\documentclass[a4paper]{article}
\usepackage{latexsym}
\usepackage{amssymb}
\begin{document}
\title{Constructing $(\omega _1 ,\beta )$-morasses for $\omega _1 \leq \beta$}
\author{Bernhard Irrgang}
\date{December 16, 2010}
\maketitle
\begin{abstract}
Let $\kappa \in Card$ and $L_\kappa [X]$ be such that the fine structure theory, condensation and $Card ^{L_\kappa [X]}=Card \cap \kappa$ hold. Then it is possible to prove the existence of morasses. In particular, I will prove that there is a $\kappa$-standard morass, a notion that I introduced in a previous paper. This shows the consistency of $(\omega _1 , \beta )$-morasses for all $\beta \geq \omega _1$.
\end{abstract}  
\section{Introduction}
R. Jensen formulated in the 1970's the concept of an $(\omega_\alpha,\beta)$-morass whereby objects of size $\omega_{\alpha+\beta}$ could be constructed by a directed system of objects of size less than $\omega_\alpha$. He defined the notion of an $(\omega_\alpha,\beta)$-morass only for the case that $\beta < \omega_\alpha$.  I introduced in a previous paper [Irr2] a definition of an $(\omega_1,\beta)$-morass for the case that $\omega_1 \leq \beta$.  
\smallskip\\
This definition of an $(\omega_1,\beta)$-morass for the case that $\omega_1 \leq \beta$ seems to be an axiomatic description of  the condensation property of G\"odel's constructible universe $L$ and the whole fine structure theory of it. I was, however, not able to formulate and prove this fact in form of a mathematical statement. Therefore, I defined a seemingly innocent strengthening of the notion of an $(\omega_1,\beta)$-morass, which I actually expect to be equivalent to the notion of $(\omega_1,\beta)$-morass. I call this strengthening an $\omega_{1+\beta}$-standard morass. As will be seen, if we construct a morass in the usual way in $L$, the properties of a standard morass hold automatically. 
\smallskip\\
Using the notion of a standard morass, I was able to prove a theorem which can be interpreted as saying that standard morasses fully cover the condensation property and fine structure of $L$. More precisely, I was able to show the following [Irr2]
\medskip\\
{\bf Theorem}
\smallskip\\
Let $\kappa \geq \omega_1$ be a cardinal and assume that a $\kappa$-standard morass exists. Then there exists a predicate $X$ such that $Card \cap \kappa=Card^{L_\kappa[X]}$ and $L_\kappa[X]$ satisfies amenability, coherence and condensation.
\medskip\\
Let me explain this. The predicate $X$ is a sequence $X=\langle X_\nu \mid \nu \in S^X\rangle$ where $S^X \subseteq Lim \cap \kappa$, and $L_\kappa[X]$ is endowed with the following hierarchy: Let $I_\nu =\langle J^X_\nu , X \upharpoonright \nu \rangle$ for $\nu \in Lim -S^X$ and $I_\nu = \langle J^X_\nu , X \upharpoonright \nu ,X_\nu \rangle$ for $\nu \in S^X$ where $X_\nu \subseteq J^X_\nu$ and
\smallskip

$J^X_0 = \emptyset$
\smallskip

$J^X_{\nu +\omega } = rud(I^X_\nu )$
\smallskip

$J^X_\lambda = \bigcup \{ J^X_\nu \mid \nu \in \lambda \}$ for $\lambda \in Lim^2:=Lim(Lim)$,
\smallskip \\
where $rud(I^X_\nu )$ is the rudimentary closure of $J^X_\nu \cup \{ J^X_\nu \}$ relative to $X \upharpoonright \nu$ if $\nu \in Lim - S^X$ and relative to $X \upharpoonright \nu$ and $X_\nu$ if $\nu \in S^X$. Now, the properties of $L_\kappa[X]$ are defined as follows:
\smallskip\\
{\bf (Amenability)} The structures $I_\nu$ are amenable.
\smallskip\\
{\bf (Coherence)} If $\nu \in S^X$, $H \prec _1 I_\nu$ and $\lambda = sup(H \cap On)$, then $\lambda \in S^X$ and $X_\lambda = X_\nu \cap J^X_\lambda$.
\smallskip \\
{\bf (Condensation)} If $\nu \in S^X$ and $H \prec_1 I_\nu$, then there is some $\mu \in S^X$ such that $H \cong I_\mu$.
\smallskip\\
Moreover, if we let $\beta (\nu )$ be the least $\beta$ such that $J^X_{\beta +\omega} \models \nu$ singular, then $S^X =\{ \beta (\nu )\mid \nu$ singular in $I_\kappa \}.$
\medskip\\
As will be seen, these properties suffice to develop the fine structure theory. In this sense, the theorem shows indeed what I claimed. In the present paper, I shall show the converse:
\medskip\\
{\bf Theorem}
\smallskip\\
If $L_\kappa [X]$, $\kappa \in Card$, satisfies condensation, coherence, amenability, $S^X =\{ \beta (\nu )\mid \nu$ singular in $I_\kappa \}$ and $Card ^{L_\kappa [X]}=Card \cap \kappa$, then there is a $\kappa$-standard morass.
\medskip\\
Since $L$ itself satisfies the properties of $L_\kappa[X]$, this also shows that the existence of $\kappa$-standard morasses and $(\omega_1,\beta)$-morasses is consistent for all $\kappa \geq\omega_2$ and all $\beta\geq \omega_1$. 
\smallskip\\
Most results that can be proved in $L$ from condensation and the fine structure theory also hold in the structures $L_\kappa[X]$ of the above form. As examples, I proved in my dissertation the following two theorems whose proofs can also be seen as applications of morasses:
\medskip\\
{\bf Theorem}
\smallskip\\
Let $\lambda \geq \omega_1$ be a cardinal, $S^X \subseteq Lim \cap \lambda$, $Card \cap \lambda = Card^{L_\lambda[X]}$ and  $X=\langle X_\nu \mid \nu \in S^X\rangle$ be a sequence such that amenability, coherence, condensation and $S^X =\{ \beta (\nu )\mid \nu$ singular in $I_\kappa \}$ hold. Then $\Box_\kappa$ holds for all infinite cardinals $\kappa < \lambda$.
\medskip\\
{\bf Theorem}
\smallskip\\
Let $S^X \subseteq Lim$ and  $X=\langle X_\nu \mid \nu \in S^X\rangle$ be a sequence such that amenability, coherence, condensation and $S^X =\{ \beta (\nu )\mid \nu$ singular in $L[X] \}$ hold.  Then the weak covering lemma holds for $L[X]$. That is, if there is no non-trival, elementary embedding $\pi:L[X] \rightarrow L[X]$, $\kappa \in Card^{L[X]}-\omega_2$ and $\tau = (\kappa^+)^{L[X]}$, then
$$\tau < \kappa^+ \quad \Rightarrow \quad cf(\tau)=card(\kappa).$$ 
\smallskip\\
The present paper is a part of my dissertation [Irr1]. I thank Dieter Donder for being my adviser, Hugh Woodin for an invitation to Berkeley, where part of the work was done, and the DFG-Graduiertenkolleg ``Sprache, Information, Logik`` in Munich for their support.

\section{The inner model $L[X]$}
We say a function $f:V^n \rightarrow V$ is rudimentary for some structure $\frak{W} = \langle W, X_i \rangle$ if it is generated by the following schemata:
\smallskip

$f(x_1, \dots , x_n) = x_i$ for $1 \leq i \leq n$
\smallskip 

$f(x_1, \dots , x_n) = \{ x_i, x_j \}$ for $1 \leq i,j \leq n$
\smallskip

$f(x_1, \dots , x_n) = x_i - x_j$ for $1 \leq i,j \leq n$ 
\smallskip

$f(x_1 , \dots , x_n) = h(g_1(x_1, \dots , x_n), \dots , g_n(x_1, \dots , x_n))$

where $h, g_1, \dots , g_n$ are rudimentary
\smallskip

$f(y, x_2, \dots , x_n) = \bigcup \{ g(z, x_2, \dots , x_n) \mid z \in y \}$

where $g$ is rudimentary
\smallskip 

$f(x_1, \dots , x_n) = X_i \cap x_j$ where $1 \leq j \leq n$.
\bigskip\\
{\bf Lemma 1}
\smallskip \\
A function is rudimentary iff it is a composition of the following functions: 
\smallskip

$F_0(x,y)=\{ x,y \}$
\smallskip

$F_1(x,y)=x-y$
\smallskip

$F_2(x,y)=x \times y$
\smallskip

$F_3(x,y)=\{ \langle u,z,v \rangle  \mid z \in x$ $and$ $\langle u,v \rangle \in y \}$
\smallskip

$F_4(x,y)=\{ \langle z,u,v \rangle  \mid z \in x$ $and$ $\langle u,v \rangle \in y \}$
\smallskip

$F_5(x,y)=\bigcup x$
\smallskip

$F_6(x,y)=dom(x)$
\smallskip

$F_7(x,y)=\in \cap (x \times x)$
\smallskip

$F_8(x,y)=\{ x[\{ z\} ] \mid z \in y \}$
\smallskip

$F_{9+i}(x,y)=x \cap X_i$ for the predicates $X_i$ of the structure under consideration.
\smallskip \\
{\bf Proof:} See, for example, in [Dev2].  $\Box$
\medskip \\
A relation $R \subseteq V^n$ is called rudimentary if there is a rudimentary function $f:V^n \rightarrow V$ such that $R(x_i) \Leftrightarrow f(x_i) \neq \emptyset$.
\medskip \\
{\bf Lemma 2}
\smallskip \\
Every relation that is $\Sigma _0$ over the considered structure is rudimentary.
\smallskip \\
{\bf Proof:} Let $\chi _R$ be the characteristic function of $R$. The claim follows from the facts (i)-(vi):
\smallskip \\
(i) $R$ rudimentary $\Leftrightarrow$ $\chi _R$ rudimentary.
\smallskip \\
$\Leftarrow$ is clear. Conversely, $\chi _R=\bigcup \{ g(y)\mid y \in f(x_i)\}$ where $g(y)=1$ is constant and $R(x_i)$ $\Leftrightarrow$ $f(x_i) \neq \emptyset$.
\smallskip \\
(ii) If $R$ is rudimentary, then $\neg R$ is also rudimentary.
\smallskip \\
Since $\chi _{\neg R}=1-\chi _R$.
\smallskip \\
(iii) $x \in y$ and $x=y$ are rudimentary.
\smallskip \\
By $x \notin y$ $\Leftrightarrow$ $\{ x\} - y\neq\emptyset$ , $x \neq y$ $\Leftrightarrow$ $(x-y) \cup (y-x) \neq \emptyset$ and (ii).
\smallskip \\
(iv) If $R(y,x_i)$ is rudimentary, then $(\exists z \in y)R(z,x_i)$ and $(\forall z \in y)R(z,x_i)$ are rudimentary. 
\smallskip \\ 
If $R(y,x_i)$ $\Leftrightarrow$ $f(y,x_i) \neq \emptyset$, then $(\exists z \in y)R(z,x_i)$ $\Leftrightarrow$ $\bigcup \{f(z,x_i) \mid z \in y \}\neq \emptyset$. The second claim follows from this by (ii).
\smallskip \\
(v) If $R_1 , R_2 \subseteq V^n$ are rudimentary, then so are $R_1 \vee R_2$ and $R_1 \wedge R_2$.
\smallskip \\
Because $f(x,y)=x \cup y$ is rudimentary, $(R_1 \vee R_2)(x_i )$ $\Leftrightarrow$ $\chi _{R_1}(x_i) \cup \chi _{R_2}(x_i) \neq \emptyset$ is rudimentary. The second claim follows from that by (ii).
\smallskip \\
(vi) $x \in X_i$ is rudimentary.
\smallskip \\
Since $\{ x\} \cap X_i \neq \emptyset$ $\Leftrightarrow$ $x \in X_i$. $\Box$
\medskip \\
For a converse of this lemma, we define:
\smallskip \\
A function $f$ is called simple if $R(f(x_i),y_k)$ is $\Sigma _0$ for every $\Sigma _0$-relation $R(z,y_k)$.
\medskip \\
{\bf Lemma 3}
\smallskip \\
A function $f$ is simple iff
\smallskip

(i) $z \in f(x_i)$ is $\Sigma _0$
\smallskip

(ii) $A(z)$ is $\Sigma _0$ $\Rightarrow$ $(\exists z \in f(x_i))A(z)$ is $\Sigma _0$.
\smallskip \\
{\bf Proof:} If $f$ is simple, then (i) and (ii) hold, because these are instances of the definition. The converse is proved by induction on $\Sigma _0$-formulas. E.g. if $R(z,y_k)$ $:\Leftrightarrow$ $z=y_k$, then $R(f(x_i),y_k)$ $\Leftrightarrow$ $f(x_i)=y_k$ $\Leftrightarrow$ $(\forall z \in f(x_i))(z \in y_k)$ $and$ $(\forall z \in y_k)(z \in f(x_i))$. Thus we need (i) and (ii). The other cases are similar. $\Box$
\medskip \\
{\bf Lemma 4}
\smallskip \\
Every rudimentary function is $\Sigma _0$ in the parameters $X_i$.
\smallskip \\
{\bf Proof:} By induction, one proves that the rudimentary functions that are generated without the schema $f(x_1,\dots ,x_n)=X_i \cap x_j$ are simple. For this, one uses lemma 3. But since the function $f(x,y)=x\cap y$ is one of those, the claim holds. $\Box$
\medskip \\
Thus every rudimentary relation is $\Sigma _0$ in the parameters $X_i$, but not necessaryly $\Sigma _0$ with the $X_i$ as predicates. An example is the relation $\{ x,y\} \in X_0$. 
\medskip \\
A structure is said to be rudimentary closed if its underlying set is closed under all rudimentary functions.
\medskip \\
{\bf Lemma 5}
\smallskip \\
If $\frak{W}$ is rudimentary closed and $H \prec _1 \frak{W}$, then $H$ and the collapse of $H$ are also rudimentary closed.
\smallskip \\
{\bf Proof:} That is clear, since the functions $F_0,\dots ,F_{9+i}$ are $\Sigma _0$ with the predicates $X_i$. $\Box$
\medskip \\
Let $T_N$ be the set of $\Sigma _0$ formulae of our language $\{ \in ,X_1, \dots , X_N\}$ having exactly one free variable. By lemma 2, there is a rudimentary function $f$ for every $\Sigma _0$ formula $\psi$ such that $\psi (x_\star )$ $\Leftrightarrow$ $f(x_\star ) \neq \emptyset$. By lemma 1, we have
\begin{tabbing}
where \= blabla \kill \\
\> $x_0=f(x_\star )=F_{k_1}(x_1,x_2)$ \\
where \> $x_1=F_{k_2}(x_3,x_4)$ \\
\> $x_2=F_{k_3}(x_5,x_6)$ \\
and \> $x_3=\dots$ \\
\end{tabbing}
Of course, $x_\star$ appears at some point.
\medskip \\
Therefore, we may define an effective G\"odel coding
$$T_N \rightarrow G,\psi _u \mapsto u$$
as follows ($m,n$ possibly $=\star$): 
$$\langle k,l,m,n \rangle \in u :\Leftrightarrow x_k = F_l(x_m,x_n).$$
Let $\models ^{\Sigma _0} _\frak{W}(u,x_\star )$ $:\Leftrightarrow$
\medskip

$\psi _u$ is a $\Sigma _0$ formula with exactly one free variable
\smallskip 

and $\frak{W} \models \psi _u (x_\star )$.
\bigskip\\ \\
{\bf Lemma 6}
\smallskip \\
If $\frak{W}$ is transitive and rudimentary closed, then $\models ^{\Sigma _0} _\frak{W}(x,y)$ is $\Sigma _1$-definable over $\frak{W}$. The definition of $\models ^{\Sigma _0} _\frak{W}(x,y)$ depends only on the number of predicates of $\frak{W}$. That is, it is uniform for all structures of the same type.
\smallskip \\
{\bf Proof:} Whether $\models ^{\Sigma _0} _\frak{W}(u,x_\star )$ holds, may be computed directly. First, one computes the $x_k$ which only depend on $x_\star$. For those $k$, $\langle k,l,\star ,\star \rangle \in u$. Then one computes the $x_i$ which only depend on $x_m$ and $x_n$ such that $m,n \in \{ k \mid \langle k,l,\star ,\star \rangle \in u \}$ -- etc. Since $\frak{W}$ is rudimentary closed, this process only breaks off, when one has computed $x_0=f(x_\star )$. And $\models ^{\Sigma _0} _\frak{W}(u,x_\star )$ holds iff $x_0=f(x_\star )\neq\emptyset$.

More formally speaking: $\models ^{\Sigma _0} _\frak{W}(u,x_\star )$ holds iff there is some sequence $\langle x_i \mid i \in d \rangle$, $d=\{ k \mid \langle k,l,m,n\rangle \in u \}$ such that
\smallskip

$\langle k,l,m,n \rangle \in u$ $\Rightarrow$ $x_k=F_l(x_m,x_n)$
\smallskip

and $x_0 \neq \emptyset$.
\smallskip \\
Hence $\models ^{\Sigma _0} _\frak{W}$ is $\Sigma _1$. $\Box$
\medskip \\
If $\frak{W}$ is a structure, then let $rud(\frak{W} )$ be the closure of $W \cup \{ W\}$ under the functions which are rudimentary for $\frak{W}$.
\medskip\\
{\bf Lemma 7}
\smallskip \\
If $\frak{W}$ is transitive, then so is $rud(\frak{W} )$.
\smallskip \\
{\bf Proof:} By induction on the definition of the rudimentary functions. $\Box$ 
\medskip \\
{\bf Lemma 8}
\smallskip \\
Let $\frak{W}$ be a transitive structure with underlying set $W$. Then 
$$rud(\frak{W} ) \cap \frak{P} (W)= Def(\frak{W} ).$$
{\bf Proof:} First, let $A \in Def(\frak{W})$. Then $A$ is $\Sigma _0$ over $\langle W \cup \{ W \} ,X_i \rangle$, i.e. there are parameters $p_i \in W \cup \{ W\}$ and some $\Sigma _0$ formula $\varphi$ such that $x\in A$ $\Leftrightarrow$ $\varphi (x,p_i)$. But by lemma 2, every $\Sigma _0$ relation is rudimentary. Thus there is a rudimentary function $f$ such that $x \in A$ $\Leftrightarrow$ $f(x,p_i) \neq \emptyset$. Let $g(z,x)=\{ x\}$ and define $h(y,x)=\bigcup\{ g(z,x) \mid z \in y\}$. Then $h(f(x,p_i),x) = \bigcup\{ g(z,x)\mid z \in f(x,p_i)\}$ is rudimentary, $h(f(x,p_i),x)=\emptyset$ if $x \notin A$ and $h(f(x,p_i),x)=\{x\}$ if $x \in A$. Finally, let $H(y,p_i)=\bigcup\{ h(f(x,p_i),x)\mid x \in y\}$. Then $H$ is rudimentary and $A=H(W,p_i)$. So we are done.

Conversely, let $A \in rud(\frak{W} ) \cap \frak{P} (W)$. Then there is a rudimentary function $f$ and some $a \in W$ such that $A=f(a,W)$. By lemma 4 and lemma 3, there exists a $\Sigma _0$ formula $\psi$ such that $x \in f(a,W)$ $\Leftrightarrow$ $\psi (x,a,W,X_i)$. By $\Sigma _0$ absoluteness, $A=\{ x \in W \mid W\cup \{W,X_i\} \models \psi (x,a,W,X_i)\}$, since $X_i \subseteq W$. Therefore, there is a formula $\varphi$ such that $A=\{ x \in W \mid \frak{W} \models \varphi (x,a)\}$. $\Box$
\medskip \\
Let $\kappa \in Card -\omega _1$, $S^X \subseteq Lim \cap \kappa$ and $\langle X_\nu \mid \nu \in S^X \rangle$ be a sequence.
\smallskip \\
For $\nu \in Lim -S^X$, let $I_\nu =\langle J^X_\nu , X \upharpoonright \nu \rangle$ and let $I_\nu = \langle J^X_\nu , X \upharpoonright \nu ,X_\nu \rangle$ for $\nu \in S^X$ such that
\smallskip \\
$X_\nu \subseteq J^X_\nu$ where
\smallskip

$J^X_0 = \emptyset$
\smallskip

$J^X_{\nu +\omega } = rud(I_\nu )$
\smallskip

$J^X_\lambda = \bigcup \{ J^X_\nu \mid \nu \in \lambda \}$ if $\lambda \in Lim^2:=Lim(Lim)$.
\medskip \\
Obviously, $L_\kappa [X]=\bigcup \{ J^X_\nu \mid \nu \in \kappa  \}$.
\smallskip \\
We say that $L_\kappa [X]$ is amenable if $I_\nu$ is rudimentary closed for all $\nu \in S^X$.
\bigskip\\
{\bf Lemma 9}
\smallskip 

(i) Every $J_\nu ^X$ is transitive
\smallskip 

(ii) $\mu < \nu$ $\Rightarrow$ $J_\mu ^X \in J_\nu ^X$
\smallskip 

(iii) $rank(J_\nu ^X)=J_\nu ^X \cap On = \nu$
\smallskip \\
{\bf Proof:} That are three easy proofs by induction. $\Box$
\medskip \\
Sometimes we need levels between $J^X_\nu$ and $J^X_{\nu +\omega}$. To make those transitive, we define
\smallskip 

$G_i(x,y,z)=F_i(x,y)$ for $i \leq 8$
\smallskip

$G_9(x,y,z)=x \cap X$
\smallskip

$G_{10}(x,y,z)=\langle x,y \rangle$
\smallskip 

$G_{11}(x,y,z)=x[y]$
\smallskip

$G_{12}(x,y,z)=\{ \langle x,y \rangle \}$
\smallskip 

$G_{13}(x,y,z)=\langle x,y,z\rangle$
\smallskip

$G_{14}(x,y,z)=\{ \langle x,y\rangle ,z\}$.
\medskip \\
Let
\smallskip

$S_0=\emptyset$
\smallskip 

$S_{\mu +1}=S_\mu \cup \{ S_\mu \} \cup \bigcup \{ G_i [(S_\mu \cup \{ S_\mu \} )^3] \mid i \in 15\}$
\smallskip

$S_\lambda =\bigcup \{ S_\mu \mid \mu \in \lambda \}$ if $\lambda \in Lim$.
\medskip \\
{\bf Lemma 10}
\smallskip \\
The sequence $\langle I_\mu \mid \mu \in Lim \cap \nu  \rangle$ is (uniformly) $\Sigma _1$-definable over $I_\nu$.
\smallskip \\
{\bf Proof:} By definition $J^X_\mu=S_\mu$ for $\mu \in Lim$, that is, the sequence $\langle J^X_\mu \mid \mu \in Lim \cap \nu  \rangle$ is the solution of the recursion defining $S_\mu$ restricted to $Lim$. Since the recursion condition is $\Sigma _0$ over $I_\nu$, the solution is $\Sigma _1$. It is $\Sigma _1$ over $I_\nu$ if the existential quantifier can be restricted to $J^X_\nu$. Hence we must prove $\langle S_\mu \mid \mu \in \tau \rangle \in J^X_\nu$ for $\tau \in \nu$. This is done by induction on $\nu$. The base case $\nu =0$ and the limit step are clear. For the successor step, note that $S_{\mu +1}$ is a rudimentary function of $S_\mu$ and $\mu$, and use the rudimentary closedness of $J^X_\nu$. $\Box$
\medskip \\
{\bf Lemma 11}
\smallskip \\
There are well-orderings $<_\nu$ of the sets $J_\nu^X$ such that
\smallskip 

(i) $\mu < \nu$ $\Rightarrow$ $<_\mu \subseteq <_\nu$
\smallskip 

(ii) $<_{\nu +1}$ is an end-extension of $<_\nu$
\smallskip 

(iii) The sequence $\langle <_\mu \mid \mu \in Lim \cap \nu  \rangle$ is (uniformly) $\Sigma _1$-definable over $I_\nu$.
\smallskip 

(iv) $<_\nu$ is (uniformly) $\Sigma _1$-definable over $I_\nu$.
\smallskip 

(v) The function $pr_\nu (x)= \{ z \mid z<_\nu x \}$ is (uniformly) $\Sigma _1$-definable over $I_\nu$.
\smallskip \\
{\bf Proof:} Define well-orderings $<_\mu$ of $S_\mu$ by recursion: \\
\begin{tabbing}
(III)\= (1) \= (a) \= \kill \\
(I)\> \> $<_0=\emptyset$ \\
\\
(II)\> (1) For $x,y\in S_\mu$, let $x<_{\mu +1}y$ $\Leftrightarrow$ $x<_\mu y$ \\
\> (2) $x \in S_\mu$ $and$ $y \notin S_\mu$ $\Rightarrow$ $x<_{\mu +1}y$ \\
\> \> $y \in S_\mu$ $and$ $x \notin S_\mu$ $\Rightarrow$ $y<_{\mu +1}x$ \\
\> (3) If $x,y \notin S_\mu$, then there is an $i \in 15$ and $x_1,x_2,x_3 \in S_\mu$ such that\\ 
\> \> $x=G_i(x_1,x_2,x_3)$. And there is a $j \in 15$ and $y_1,y_2,y_3 \in S_\mu$\\ 
\> \> such that $y=G_j(y_1,y_2,y_3)$. First, choose $i$ and $j$ minimal, then\\
\> \> $x_1$ and $y_1$, then $x_2$ and $y_2$, and finally $x_3$ and $y_3$. \\
\> \> Set: \\
\> \> (a) $x<_{\mu +1} y$ if $i<j$ \\
\> \> \> $y<_{\mu +1} x$ if $j<i$ \\
\> \> (b) $x<_{\mu +1} y$ if $i=j$ $and$ $x_1<_\mu y_1$ \\
\> \> \> $y<_{\mu +1} x$ if $i=j$ $and$ $y_1<_\mu x_1$ \\
\> \> (c) $x<_{\mu +1} y$ if $i=j$ $and$ $x_1= y_1$ $and$ $x_2<_\mu y_2$ \\
\> \> \> $y<_{\mu +1} x$ if $i=j$ $and$ $x_1=y_1$ $and$ $y_2<_\mu x_2$ \\
\> \> (d) $x<_{\mu +1} y$ if $i=j$ $and$ $x_1= y_1$ $and$ $x_2= y_2$ $and$ $x_3 <_\mu y_3$ \\
\> \> \> $y<_{\mu +1} x$ if $i=j$ $and$ $x_1=y_1$ $and$ $y_2= x_2$ $and$ $y_3 <_\mu x_3$\\
\\
(III) \> \> $<_\lambda =\bigcup \{ <_\mu \mid \mu \in \lambda \}$ \\
\end{tabbing}
The properties (i) to (v) are obvious. For the $\Sigma _1$-definability, one needs the argument from lemma 10. $\Box$
\medskip \\
{\bf Lemma 12}
\smallskip \\
The rudimentary closed $\langle J^X_\nu ,X \upharpoonright \nu , A \rangle$ have a canonical $\Sigma _1$-Skolem function $h$.
\smallskip \\
{\bf Proof:} Let $\langle \psi _i \mid i \in \omega \rangle$ be an effective enumeration of the $\Sigma _0$ formulae with three free variables. Intuitively, we would define:
$$h(i,x) \quad \simeq \quad (z)_0$$
for
$$\hbox{the} \quad <_\nu \hbox{-least} \quad z \in J_\nu ^X \quad \hbox{such that} \quad \langle J^X_\nu ,X \upharpoonright \nu ,A \rangle \models \psi _i ((z)_0 , x, (z)_1).$$
Formally, we define:

By lemma 11 (v), let $\theta$ be a $\Sigma _0$ formula such that
$$w=\{ v \mid v <_\nu z\} \quad \Leftrightarrow \quad \langle J^X_\nu ,X \upharpoonright \nu , A \rangle \models (\exists t) \theta (w,z,t).$$
Let $u_i$ be the G\"odel coding of
$$\theta ((s)_1,(s)_0,(s)_2)$$ $$ \wedge \quad \psi _i(((s)_0)_0,(s)_3,((s)_0)_1) \quad \wedge \quad (\forall v \in (s)_1) \neg \psi _i((v)_0,(s)_3,(v)_1)$$
and
\smallskip 

$y=h(i,x) \quad \Leftrightarrow$
$$(\exists s)(((s_0)_0=y \quad \wedge \quad (s)_3=x \quad \wedge \quad \models ^{\Sigma _0}_{\langle J^X_\nu ,X\upharpoonright \nu  ,A \rangle} (u_i,s)).$$
This has the desired properties. Note lemma 6! $\Box$
\medskip \\
I will denote this $\Sigma _1$-Skolem function by $h_{\nu ,A}$. Let $h_\nu := h_{\nu ,\emptyset}$.
\medskip \\
Let us say that $L_\kappa [X]$ has condensation if the following holds:
\smallskip\\
If $\nu \in S^X$ and $H \prec_1 I_\nu$, then there is some $\mu \in S^X$ such that $H \cong I_\mu$.
\medskip \\
From now on, suppose that $L_\kappa [X]$ is amenable and has condensation.
\medskip\\
Set $I^0_\nu =\langle J^X_\nu ,X \upharpoonright \nu \rangle$ for all $\nu \in Lim \cap \kappa$.
\medskip \\
{\bf Lemma 13} (G\"odel's pairing function)
\smallskip \\
There is a bijection $\Phi :On^2 \rightarrow On$ such that $\Phi (\alpha ,\beta ) \geq \alpha , \beta$ for all $\alpha ,\beta$ and $\Phi ^{-1} \upharpoonright \alpha$ is uniformly $\Sigma _1$-definable over $I^0_\alpha$ for all $\alpha \in Lim$.
\smallskip \\
{\bf Proof:} Define a well-ordering $<^\star$ on $On^2$ by

$\langle \alpha , \beta \rangle <^\star \langle \gamma , \delta \rangle$ \\
iff

$max(\alpha , \beta ) < max(\gamma , \delta )$ or

$max(\alpha , \beta ) = max(\gamma , \delta )$ $and$ $\alpha < \gamma$ or

$max(\alpha , \beta ) = max(\gamma , \delta )$ $and$ $\alpha = \gamma$ $and$ $\beta < \delta$.
\\
Let $\Phi :\langle On^2, <^\star \rangle \cong\langle  On,<\rangle$. Then $\Phi$ may be defined  by the recursion

$\Phi (0, \beta ) = sup \{ \Phi (\nu , \nu ) \mid \nu < \beta \}$

$\Phi (\alpha , \beta ) = \Phi (0,\beta ) + \alpha$ if $\alpha < \beta$

$\Phi (\alpha , \beta ) = \Phi (0,\alpha ) + \alpha +\beta$ if $\alpha \geq \beta$.
\\
$\Box$
\smallskip \\
So there is a uniform map from $\alpha$ onto $\alpha \times \alpha$ for all $\alpha$ that are closed under G\"odel's pairing function. Such a map exists for all $\alpha \in Lim$. But then we have to give up uniformity.
\medskip \\
{\bf Lemma 14}
\smallskip \\
For all $\alpha \in Lim$, there exists a function from $\alpha$ onto $\alpha \times \alpha$ that is $\Sigma _1$-definable over $I^0_\alpha$.
\smallskip \\
{\bf Proof} by induction on $\alpha \in Lim$. If $\alpha$ is closed under G\"odel's pairing fuction, then lemma 13 does the job. Therefore, if $\alpha = \beta + \omega$ for some $\beta \in Lim$, we may assume $\beta \neq 0$. But then there is some over $I^0_\alpha$ $\Sigma _1$-definable bijection $j:\alpha \rightarrow \beta$. And by the induction hypothesis, there is an over $I^0_\beta$ $\Sigma _1$-definable function from $\beta$ onto $\beta \times \beta$. Thus there exists a $\Sigma _1$ formula $\varphi (x,y,p)$ and a parameter $p \in J^X_\beta$ such that there is some $x \in \beta$ satisfying $\varphi (x,y,p)$ for all $y \in \beta \times \beta$. So we get an over $I^0_\beta$ $\Sigma _1$-definable injective function $g: \beta \times \beta \rightarrow \beta$ from the $\Sigma _1$-Skolem function. Hence $f(\langle \nu , \tau \rangle ) = g( \langle j( \nu ), j(\tau ) \rangle )$ defines an injective function $f: \alpha ^2 \rightarrow \beta$ which is $\Sigma _1$-definable over $I^0_\alpha$. An $h$ which is as needed may be defined by 

$h(\nu ) = f^{-1} (\nu )$ if $\nu \in rng(f)$

$h(\nu ) = \langle 0,0 \rangle$ else.
\\
For $rng(f)=rng(g) \in J^X_\alpha$.

Now, assume $\alpha \in Lim^2$ is not closed under G\"odel's pairing function. Then $\nu , \tau \in \alpha$ for $\langle \nu , \tau \rangle = \Phi ^{-1}(\alpha )$, and $c:=\{ z \mid z <^\star \langle \nu , \tau \rangle \}$ lies in $J^X_\alpha$. Thus $\Phi ^{-1} \upharpoonright c:c \rightarrow \alpha$ is an over $I^0_\alpha$ $\Sigma _1$-definable bijection. Pick a $\gamma \in Lim$ such that $\nu ,\tau < \gamma$. Then $\Phi ^{-1} \upharpoonright \alpha : \alpha \rightarrow \gamma ^2$ is an over $I^0_\alpha$ $\Sigma _1$-definable injective function. Like in the first case, there exists an injective function $g:\gamma \times \gamma \rightarrow \gamma$ in $J^X_\alpha$ by the induction hypothesis. So $f(\langle \xi , \zeta \rangle ) = g(\langle g \Phi ^{-1}(\xi ),g\Phi ^{-1}(\zeta )) \rangle )$ defines an over $I^0_\alpha$ $\Sigma _1$-definable bijection $f:\alpha ^2 \rightarrow d$ such that $d:=g[g[c] \times g[c]]$. Again, we define $h$ by

$h(\xi )=f^{-1}(\xi )$ if $\xi \in d$

$h(\xi )=\langle 0,0 \rangle$ else. $\Box$   
\medskip \\
{\bf Lemma 15}
\smallskip \\
Let $\alpha \in Lim - \omega +1$. Then there is some over $I^0_\alpha$ $\Sigma _1$-definable function from $\alpha$ onto $J^X_\alpha$. This function is uniformly definable for all $\alpha$ closed under G\"odel's pairing function.
\smallskip \\
{\bf Proof:} Let $f:\alpha \rightarrow \alpha \times \alpha$ be a surjective function which is $\Sigma _1$-definable over $I^0_\alpha$ with parameter $p$. Let $p$ be minimal with respect to the canonical well-ordering such that such an $f$ exists. Define$f^0,f^1$ by $f(\nu )=\langle f^0(\nu ),f^1(\nu )\rangle$ and, by induction, define $f_1=id \upharpoonright \alpha$ and $f_{n+1} (\nu )= \langle f^0(\nu ),f_n \circ f^1(\nu ) \rangle$. Let $h:=h_\alpha$ be the canonical $\Sigma _1$-Skolem function and $H=h[\omega \times (\alpha \times \{ p \} )]$. Then $H$ is closed under ordered pairs. For, if $y_1=h(j_1,\langle \nu _1,p\rangle )$, $y_2=h(j_2,\langle \nu _2,p\rangle )$ and $\langle \nu _1,\nu _2\rangle =f(\tau )$, then $\langle y_1,y_2\rangle$ is $\Sigma _1$-definable over $I^0_\alpha$ with the parameters $\tau ,p$. Hence it is in $H$. Since $H$ is closed under ordered pairs, we have $H \prec _1 I^0_\alpha$. Let $\sigma : H \rightarrow I^0_\beta$ be the collapse of $H$. Then $\alpha = \beta$, because $\alpha \subseteq H$ and $\sigma \upharpoonright \alpha = id \upharpoonright \alpha$. Thus $\sigma [f]=f$, and $\sigma [f]$ is $\Sigma _1$-definable over $I^0_\alpha$ with the parameter $\sigma (p)$. Since $\sigma$ is a collapse, $\sigma (p) \leq p$. So $\sigma (p) =p$ by the minimality of $p$. In general, $\pi (h(i,x)) \simeq h(i,\pi (x))$ for $\Sigma _1$-elementary $\pi$. Therefore, $\sigma (h(i,\langle \nu ,p \rangle)) \simeq h(i,\langle \nu ,p\rangle)$ holds in our case for all $i \in \omega$ and $\nu \in \alpha$. But then $\sigma \upharpoonright  H = id \upharpoonright H$ and $H=J^X_\alpha$. Thus we may define the needed surjective map by $g \circ f_3$ where

$g(i,\nu ,\tau ) = y$ if $(\exists z\in S_\tau ) \varphi (z,y,i,\langle \nu ,p\rangle )$

$g(i,\nu ,\tau)=\emptyset$ else.
\\
Here, $S_\tau$ shall be defined as in lemma 10 and $y=h(i,x) \Leftrightarrow (\exists t \in J^X_\alpha ) \varphi (t,i,x,y)$. $\Box$
\medskip \\
Let $\langle I^0_\nu ,A \rangle :=\langle J^X_\nu ,X \upharpoonright \nu ,A \rangle$.
\smallskip \\
The idea of the fine structure theory is to code $\Sigma _n$ predicates over large structures in $\Sigma _1$ predicates over smaller structures. In the simplest case, one codes the $\Sigma _1$ information of the given structure $I^0_\beta$ in a rudimentary closed structure $\langle I^0_\rho , A \rangle$. I.e. we want to have something like:
\medskip \\
Over $I^0_\beta$, there exists a $\Sigma _1$ function $f$ such that
$$f[J_\rho ^X]=J_\beta ^X .$$
For the $\Sigma _1$ formulae $\varphi _i$,
$$\langle i,x\rangle \in A \quad \Leftrightarrow \quad I^0_\beta \models \varphi _i (f(x))$$
holds. And
$$\langle I^0_\rho , A \rangle \quad \hbox{is rudimentary closed}.$$
Now, suppose we have such an $\langle I^0_\rho , A \rangle$. Then every $B \subseteq J_\rho ^X$ that is $\Sigma _1$-definable over $I^0_\beta$ is of the form 
$$B=\{ x \mid A(i,\langle x,p \rangle ) \} \quad \hbox{for some} \quad i \in \omega , p \in J^X_\rho.$$
So $\langle I^0_\rho , B \rangle$ is rudimentary closed for all $B \in \Sigma _1(I^0_\beta ) \cap \frak{P}(J _\rho ^X)$.
\smallskip \\
The $\rho$ is uniquely determined.
\medskip \\
{\bf Lemma 16}
\smallskip \\
Let $\beta > \omega$ and $\langle I^0_\beta , B \rangle$ be rudimentary closed. Then there is at most one $\rho \in Lim$ such that
\smallskip

$\langle I^0 _\rho , C \rangle$ is rudimentary closed for all $C \in \Sigma _1(\langle I^0_\beta ,B \rangle ) \cap \frak{P}(J _\rho ^X)$
\smallskip \\
and \smallskip 

there is an over $\langle I^0_\beta , B \rangle$ $\Sigma _1$-definable function $f$ such that $f[J_\rho ^X]=J_\beta ^X$.
\smallskip \\
{\bf Proof:} Assume $\rho < \bar \rho$ both had these properties. Let $f$ be an over $\langle I^0_\beta ,B \rangle$ $\Sigma _1$-definable function such that $f[J_\rho ^X]=J_\beta ^X$ and $C=\{ x \in J_\rho ^X \mid x \not\in f(x) \}$. Then $C \subseteq J_\rho ^X$ is $\Sigma _1$-definable over $\langle I^0_\beta ,B \rangle$. So $\langle I^0_{\bar \rho},C\rangle$ is rudimentary closed. But then $C=C\cap J^X_\rho \in J^X_{\bar \rho}$. Hence there is an $x \in J^X_\rho$ such that $C=f(x)$. From this, the contradiction $x \in f(x)$ $\Leftrightarrow$ $x \in C$ $\Leftrightarrow$ $x \not\in f(x)$ follows. $\Box$ 
\medskip \\
The uniquely determined $\rho$ from lemma 16 is called the projectum of $\langle I^0_\beta , B\rangle$.
\smallskip \\
If there is some over $\langle I^0_\beta , B \rangle$ $\Sigma _1$-definable function $f$ such that $f[J_\rho ^X]=J^X_\beta$, then $h_{\beta , B}[\omega \times (J^X_\rho \times \{ p \} )] = J^X_\beta$ for a $p \in J^X_\beta$. Using the canonical function $h_{\beta ,B}$, we can define a canonical $A$:
\smallskip \\
Let $p$ be minimal with respect to the canonical well-ordering such that the above property holds. Define 
$$A=\{ \langle i,x \rangle \mid i \in \omega \quad and \quad x \in J^X_\rho \quad and \quad \langle I^0_\beta , B \rangle \models \varphi _i (x,p)\}.$$
We say $p$ is the standard parameter of $\langle I^0_\beta , B \rangle$ and $A$ the standard code of it.
\medskip \\
{\bf Lemma 17}
\smallskip \\
Let $\beta > 0$ and $\langle I^0_\beta , B\rangle$ be rudimentary closed. Let $\rho$ be the projectum and $A$ the standard code of it. Then for all $m \geq 1$, the following holds:
$$\Sigma _{1+m}(\langle I^0_\beta ,B\rangle ) \cap \frak{P}(J^X_\rho )=\Sigma _m(\langle I^0_\rho , A\rangle ).$$
{\bf Proof:} First, let $R \in \Sigma _{1+m}(\langle I^0_\beta ,B\rangle ) \cap \frak{P}(J^X_\rho )$ and let $m$ be even. Let $P$ be a relation being $\Sigma _1$-definable over $\langle I^0_\beta , B \rangle$ with parameter $q_1$ such that, for $x \in J^X_\rho$, $R(x)$ holds iff $\exists y_0 \forall y_1 \exists y_3 \dots \forall y_{m-1} P(y_i,x)$. Let $f$ be some over $\langle I^0_\beta , B \rangle$ with parameter $q_2$ $\Sigma _1$-definable function such that $f[J^X_\rho ]=J^X_\beta$. Define $Q(z_i,x)$ by $z_i,x \in J^X_\rho$ $and$ $(\exists y_i)(y_i=f(z_i)$ $and$ $P(y_i,x))$. Let $p$ be the standard parameter of $\langle I^0_\beta ,B \rangle$. Then, by definition, there is some $u \in J^X_\rho$ such that $\langle q_1,q_2 \rangle$ is $\Sigma _1$-definable in $\langle I^0_\beta , B \rangle$ with the parameters $u,p$. I.e. there is some $i \in \omega$ such that $Q(z_i,x)$ holds iff $z_i,x \in J^X_\rho$ $and$ $\langle I^0_\beta ,B \rangle \models \varphi _i(\langle z_i,x,u\rangle ,p)$ -- i.e. iff $z_i,x \in J^X_\rho$ $and$ $A(i,\langle z_i,x,u \rangle )$. Analogously there is a $j \in \omega$ and a $v \in J^X_\rho$ such that $z \in dom(f) \cap J^X_\rho$ iff $z \in J^X_\rho$ $and$ $A(j,\langle z,v \rangle )$. Abbreviate this by $D(z)$. But then, for $x \in J^X_\rho$, $R(x)$ holds iff $\exists y_0 \forall y_1 \exists y_3 \dots \forall y_{m-1}(D(z_0) \wedge D(z_2) \wedge \dots \wedge D(z_{m-2})$ $and$ $(D(z_1) \wedge D(z_3) \wedge \dots \wedge D(z_{m-1}) \Rightarrow Q(z_i,x)))$. So the claim holds. If $m$ is odd, then we proceed correspondingly. Thus $\Sigma _{1+m}(\langle I^0_\beta ,B\rangle ) \cap \frak{P}(J^X_\rho )\subseteq \Sigma _m(\langle I^0_\rho , A\rangle )$ is proved.

Conversely, let $\varphi$ be a $\Sigma _0$ formula and $q \in J^X_\rho$ such that, for all $x \in J^X_\rho$, $R(x)$ holds iff $\langle I^0_\rho ,A \rangle \models \varphi (x,q)$. Since $\langle I^0_\rho , A \rangle$ is rudimentary closed, $R(x)$ holds iff $(\exists u \in J^X_\rho ) (\exists a \in J^X_\rho )(u$ transitive $and$ $x \in u$ $and$ $q \in u$ $and$ $a=A \cap u$ $and$ $\langle u,a \rangle \models \varphi (x,q))$. Write $a=A \cap u$ as formula: $(\forall v \in a)(v \in u$ $and$ $v \in A)$ $and$ $(\forall v \in u)(v \in A \Rightarrow v \in a)$. If $m=1$, we are done provided we can show that this is $\Sigma _2$ over $\langle I^0_\beta , B \rangle$. If $m>1$, the claim follows immediately by induction.  The second part is $\Pi _1$. So we only have to prove that the first part is $\Sigma _2$ over $\langle I^0_\beta ,B \rangle$. By the definition of $A$, $v \in A$ is $\Sigma _1$-definable over $\langle I^0_\beta , B \rangle$. I.e. there is some $\Sigma _0$ formula $\psi$ and some parameter $p$ such that $v \in A \Leftrightarrow \langle I^0_\beta , B \rangle \models (\exists y)\psi (v,y,p)$. Now, we have two cases.

In the first case, there is no over $\langle I^0_\beta , B \rangle$ $\Sigma _1$-definable function from some $\gamma < \rho$ cofinal in $\beta$. Then $(\forall v \in a)(v \in A)$ is $\Sigma _2$ over $\langle I^0_\beta , B \rangle$, because some kind of replacement axiom holds, and $(\forall v \in a) (\exists y)\psi (v,y,p)$  is over $\langle I^0_\beta , B \rangle $ equivalent to $(\exists z)(\forall v \in a)(\exists y \in z)\psi (v,y,p)$. For $\rho = \omega$, this is obvious. If $\rho \neq \omega$, then $\rho \in Lim^2$ and we can pick a $\gamma < \rho$ such that $a \in J^X_\gamma$. Let $j:\gamma \rightarrow J^X_\gamma$ an over $I_\gamma$ $\Sigma _1$-definable surjection, and $g$ an over $ \langle I^0_\beta , B \rangle$ $\Sigma _1$-definable function that maps $v \in J^X_\beta$ to $g(v) \in J^X_\beta$ such that $\psi (v,g(v),p)$ if such an element exists. We can find such a function with the help of the $\Sigma _1$-Skolem function. Now, define a function $f:\gamma \rightarrow \beta$ by

$f(\nu )=$ the least $\tau < \beta$ such that $g \circ j(\nu ) \in S_\tau$ if $j(\nu ) \in a$

$f(\nu )=0$ else. 
\\
Since $f$ is $\Sigma _1$, there is, in the given case, a $\delta < \beta$ such that $f[ \gamma ] \subseteq \delta$. So we have as collecting set $z=S_\delta$, and the equivalence is clear.

Now, let us come to the second case. Let $\gamma < \rho$ be minimal such that there is some over $\langle I^0_\beta , B \rangle$ $\Sigma _1$-definable function $g$ from $\gamma$ cofinal in $\beta$. Then $(\forall v \in a)(\exists y) \psi (v,y,p)$ is equivalent to $(\forall v \in a)(\exists \nu \in \gamma )(\exists y \in S_{g(\nu )}) \psi (v,y,p)$. If we define a predicate $C \subseteq J^X_\rho$ by $\langle v,\nu \rangle \in C \Leftrightarrow y \in S_{g(\nu )}$ $and$ $\psi (v,y,p)$, then $\langle I^0_\beta , B \rangle \models (\forall v \in a)(\exists y) \psi (v,y,p)$ is equivalent to $\langle I^0_\rho , C \rangle \models (\forall v \in a)(\exists \nu \in \gamma )(\exists y)(\langle v,\nu \rangle \in C )$. But this holds iff $\langle I^0_\rho , C \rangle \models (\exists w)(w$ transitive $and$ $a,\gamma \in w$ $and$ $\langle w,C \cap w \rangle \models (\forall v \in a)(\exists \nu \in \gamma )(\exists y)(\langle v,\nu \rangle \in C \cap w)$. Since $C$ is $\Sigma _1$ over $\langle I^0_\beta ,B \rangle$, $\langle I^0_\rho , C \rangle$ is rudimentary closed by the definition of the projectum. I.e. the statement is equivalent to $\langle I^0_\rho , C \rangle \models (\exists w)(\exists c)(w$ transitive $and$ $a,\gamma \in w$ $and$ $c=C \cap w$ $and$ $\langle w,c \rangle \models (\forall v \in a)(\exists \nu \in \gamma )(\exists y)(\langle v,\nu \rangle \in c )$. So, to prove that this is $\Sigma _2$, it suffices to show that $c=C \cap w$ is $\Sigma _2$. In its full form, this is $(\forall z)(z \in a \Leftrightarrow z \in w$ $and$ $z \in C)$. But $z \in C$ is even $\Delta _1$ over $\langle I^0_\beta , B \rangle$ by the definition. So we are finished. $\Box$
\medskip\\
{\bf Lemma 18}
\smallskip \\
(a) Let $\pi : \langle J^X_{\bar \beta},X \upharpoonright \bar \beta , \bar B \rangle \rightarrow \langle J^X_\beta , X \upharpoonright \beta , B \rangle$ be $\Sigma _0$-elementary and $\pi [\bar \beta ]$ be cofinal in $\beta$. Then $\pi$ is even $\Sigma _1$-elementary.
\smallskip \\
 (b) Let $\langle J^X_{\bar \nu}, X \upharpoonright \bar \nu , \bar A \rangle$ be rudimentary closed and $\pi : \langle J^X_{\bar \nu},X \upharpoonright \bar \nu \rangle \rightarrow \langle J^Y _\nu ,Y \upharpoonright \nu \rangle$ be $\Sigma _0$-elementary and cofinal. Then there is a uniquely determined $A \subseteq J^Y_\nu$ such that $\pi : \langle J^X_{\bar \nu}, X \upharpoonright \bar \nu , \bar A \rangle \rightarrow \langle J^Y _\nu ,X \upharpoonright \nu , A \rangle$ is $\Sigma _0$-elementary and $\langle J^Y _\nu ,X \upharpoonright \nu , A \rangle$ is rudimentary closed.
\smallskip \\
{\bf Proof:} (a) Let $\varphi$ be a $\Sigma _0$ formula such that $\langle J^X_\beta , X \upharpoonright \beta , B \rangle \models (\exists z) \varphi (z, \pi (x_i))$. Since $\pi [\bar \beta ]$ is cofinal in $\beta$, there is a $\nu \in \bar\beta$ such that $\langle J^X_\beta , X \upharpoonright \beta ,B \rangle \models (\exists z \in S_{\pi (\nu )})\varphi (z, \pi (x_i))$. Here, the $S_\nu$ is defined as in lemma 10. If $\pi (S_\nu )=S_{\pi (\nu )}$, then $\langle J^X_\beta ,X \upharpoonright \beta , B \rangle \models (\exists z \in \pi (S_\nu ))\varphi (z, \pi (x_i))$. So, by the $\Sigma _0$-elementarity of $\pi$, $\langle J^X_{\bar \beta} , X \upharpoonright \bar \beta , \bar B \rangle \models (\exists z \in S_\nu )\varphi (z,x_i)$. I.e. $\langle J^X_{\bar \beta} ,X \upharpoonright \bar \beta ,\bar B \rangle \models (\exists z)\varphi (z,x_i)$. The converse is trivial.

It remains to prove $\pi (S_\nu )=S_{\pi (\nu )}$. This is done by induction on $\nu$. If $\nu =0$ or $\nu \notin Lim$, then the claim is obvious by the definition of $S_\nu$ and the induction hypothesis. So let $\lambda \in Lim$ and $M:= \pi (S_\lambda )$. Then $M$ is transitive by the $\Sigma _0$-elementarity of $\pi$. And since $\lambda \in Lim$ (i.e. $S_\lambda = J^X_\lambda$), $\langle S_\nu \mid \nu < \lambda \rangle$ is definable over $\langle J^X_\lambda , X \upharpoonright \lambda \rangle$ by (the proof of) lemma 10. Let $\varphi$ be the formula $(\forall x)(\exists \nu )(x \in S_\nu )$. Since $\pi$ is $\Sigma _0$-elementary, $\pi \upharpoonright S_\lambda : \langle J^X_\lambda , X \upharpoonright \lambda \rangle \rightarrow \langle M, (X \upharpoonright \lambda ) \cap M \rangle$ is elementary. Thus, if $\langle J^X_\lambda , X \upharpoonright \lambda \rangle \models \varphi$, then also $\langle M, (X \upharpoonright \lambda ) \cap M \rangle \models \varphi$. Since $M$ is transitive, we get $M=S_\tau$ for a $\tau \in Lim$. And, by $\pi (\lambda )=\pi (S_\lambda \cap On) = S_\tau \cap On = \tau$, it follows that $\pi (S_\lambda )=S_{\pi ( \lambda )}$.
\smallskip \\
(b) Since $\langle J^X_{\bar \nu}, X \upharpoonright \bar \nu , \bar A \rangle$ is rudimentary closed, $\bar A \cap S_\mu \in J^X_{\bar \nu}$ for all  $\mu < \bar \nu$ where $S_\mu$ is defined as in lemma 10. As in the proof of (a), $\pi (S_\mu ) = S_{\pi (\mu )}$. So we need $\pi (\bar A \cap S_\mu )=A \cap S_{\pi ( \mu )}$ to get that $\pi : \langle J^X_{\bar \nu}, X \upharpoonright \bar \nu , \bar A \rangle \rightarrow \langle J^Y _\nu ,X \upharpoonright \nu , A \rangle$ is $\Sigma _0$-elementary. Since $\pi$ is cofinal, we necessarily obtain $A=\bigcup \{ \pi (\bar A \cap S_\mu ) \mid \mu < \bar \nu \}$. But then $\langle J^Y _\nu ,X \upharpoonright \nu , A \rangle$ is rudimentary closed. For, if $x \in J^X_\nu$, we can choose some $\mu < \bar \nu$ such that $x \in S_{\pi ( \mu )}$. And $x \cap A = x \cap (A \cap S_{\pi (\mu )})=x \cap \pi (\bar A \cap S_\mu ) \in J^X_\nu$. Now, let  $\langle J^X_{\bar \nu}, X \upharpoonright \bar \nu , \bar A \rangle\models \varphi (x_i)$ where $\varphi$ is a $\Sigma _0$ formula and $u \in J^X_{\bar \nu}$ is transitive such that $x_i \in u$. Then $\langle u, X \upharpoonright \bar \nu \cap u, A \cap u \rangle \models \varphi (x_i)$ holds. Since $\pi : \langle J^X_{\bar \nu},X \upharpoonright \bar \nu \rangle \rightarrow \langle J^Y _\nu ,Y \upharpoonright \nu \rangle$ is $\Sigma _0$-elementary, $\langle \pi (u), Y \upharpoonright \nu \cap \pi (u), A \cap \pi (u) \rangle \models \varphi (\pi (x_i))$. Because $\pi (u)$ is transitive, we get $\langle J^Y _\nu ,X \upharpoonright \nu , A \rangle\models \varphi (\pi (x_i))$. This argument works as well for the converse. $\Box$ 
\medskip \\
Write $Cond_B(I^0_\beta )$ if there exists for all $H \prec _1 \langle I^0_\beta ,B \rangle$ some $\bar \beta$ and some $\bar B$ such that $H \cong \langle I^0_{\bar \beta},\bar B \rangle$.
\medskip \\
{\bf Lemma 19} (Extension of embeddings)
\smallskip \\
Let $\beta > \omega$, $m \geq 0$ and $\langle I^0_\beta , B \rangle$ be a rudimentary closed structure. Let $Cond_B(I^0_\beta )$ hold. Let $\rho$ be the projectum of $\langle I^0_\beta , B \rangle$, $A$ the standard code and $p$ the standard parameter of $\langle I^0_\beta , B \rangle$. Then $Cond _A(I^0_\rho )$ holds. And if $\langle I^0_{\bar \rho} , \bar A \rangle$ is rudimentary closed and $\pi : \langle I^0_{\bar \rho}, \bar A \rangle \rightarrow \langle I^0_\rho ,A \rangle$ is $\Sigma _m$-elementary, then there is an uniquely determined $\Sigma _{m+1}$-elementary extension $\tilde \pi : \langle I^0_{\bar \beta},\bar B \rangle \rightarrow \langle I^0_\beta , B \rangle$ of $\pi$ where $\bar \rho$ is the projectum of $\langle I^0_{\bar \beta}, \bar B \rangle$, $\bar A$ is the standard code and $\tilde \pi ^{-1}(p)$ is the standard parameter of $\langle I^0_{\bar \beta}, \bar B \rangle$.
\smallskip \\
{\bf Proof:} Let $H=h_{\beta ,B}[\omega \times (rng(\pi ) \times \{ p\})]$ $\prec _1$ $\langle I^0_\beta , B \rangle$ and $\tilde \pi : \langle I^0_{\bar \beta},\bar B \rangle \rightarrow \langle I^0_\beta , B \rangle$ be the uncollapse of $H$.
\smallskip \\
(1) $\tilde \pi$ is an extension of $\pi$
\smallskip \\
Let $\tilde \rho = sup(\pi [\bar\rho ])$ and $\tilde A=A \cap J^X_{\tilde \rho}$. Then $\pi : \langle J^X_{\bar \rho}, X \upharpoonright \bar \rho , \bar A \rangle \rightarrow \langle J^X_{\tilde \rho}, X \upharpoonright \tilde \rho , \tilde A \rangle$ is $\Sigma _0$-elementary, and by lemma 18, it is even $\Sigma _1$-elementary. We have $rng(\pi ) = H \cap J^X_{\tilde \rho}$. Obviously $rng(\pi ) \subseteq H \cap J^X_{\tilde \rho}$. So let $y \in H \cap J^X_{\tilde \rho}$. Then there is an $i \in \omega$ and an $x \in rng(\pi )$ such that $y$ is the unique $y \in J^X_\beta$ that satisfies $\langle I^0_\beta , B \rangle \models\varphi _i( \langle y,x \rangle , p)$. So by definition of $A$, $y$ is the unique $y \in J^X_{\bar \beta}$ such that $\tilde A (i,\langle y,x \rangle )$. But $x \in rng(\pi )$ and $\pi : \langle J^X_{\bar \rho}, X \upharpoonright \bar \rho , \bar A \rangle \rightarrow \langle J^X_{\tilde \rho}, X \upharpoonright \tilde \rho , \tilde A \rangle$ is $\Sigma _1$-elementary. Therefore $y \in rng(\pi )$. So we have proved that $H$ is an $\in$-end-extension of $rng(\pi )$. Since $\pi$ is the collapse of $rng(\pi )$ and $\tilde \pi$ the collapse of $H$, we obtain $\pi \subseteq \tilde \pi$.
\smallskip \\
(2) $\tilde \pi : \langle I^0_{\bar \beta},\bar B \rangle \rightarrow \langle I^0_\beta , B \rangle$  is $\Sigma _{m+1}$-elementary
\smallskip \\
We must prove $H \prec _{m+1} \langle I^0_\beta , B \rangle$. If $m=0$, this is clear. So let $m>0$ and let $y$ be $\Sigma _{m+1}$-definable in $\langle I^0_\beta , B \rangle$ with parameters from $rng(\pi ) \cup \{ p \}$. Then we have to show $y \in H$. Let $\varphi$ be a $\Sigma _{m+1}$ formula and $x_i \in rng(\pi )$ such that $y$ is uniquely determined by $\langle I^0_\beta , B \rangle \models \varphi (y,x_i,p)$. Let $\tilde h (\langle i,x \rangle ) \simeq h(i,\langle x,p \rangle )$. Then $\tilde h [J^X_\rho ]=J^X_\beta$ by the definition of $p$. So there is a $z \in J^X_\rho$ such that $y=\tilde h(z)$. If such a $z$ lies in $J^X_\rho \cap H$, then also $y \in H$, since $z,p \in H \prec _1 \langle I^0_\beta , B \rangle$. Let $D=dom(\tilde h) \cap J^X_\rho$. Then it suffices to show 
$$(\star ) \quad (\exists z_0 \in D)(\forall z_1 \in D) \dots \langle I^0_\beta , B \rangle \models \psi (\tilde h (z_i),\tilde h(z),x_i, p)$$ 
for some $z \in H\cap J^X_\rho$ where $\psi$ is $\Sigma _1$ for even $m$ and $\Pi _1$ for odd $m$ such that $\varphi (y,x_i,p) \Leftrightarrow \langle I^0_\beta ,B \rangle \models (\exists z_0)(\forall z_1) \dots \psi (z_i,y,x_i,p)$. First, let $m$ be even. Since $A$ is the standard code, there is an $i_0 \in \omega$ such that $z \in D \Leftrightarrow A(i_0,x)$ holds for all $z \in J^X_\rho$ -- and a $j_0 \in \omega$ such that, for all $z_i,z \in D$, $\langle I^0_\beta , B \rangle \models \psi (\tilde h (z_i),\tilde h(z),x_i, p)$ iff $A(j_0,\langle z_i,z,x_i \rangle )$. Thus $(\star )$ is, for $z \in J^X_\rho$, equivalent with an obvious $\Sigma _m$ formula. If $m$ is odd, then write in $(\star ) \quad \dots \neg \langle I^0_\beta , B \rangle \models \neg \psi (\dots )$. Then $\neg \psi$ is $\Sigma _1$ and we can proceed as above. Eventually $\pi : \langle I^0_{\bar \rho}, \bar A \rangle \rightarrow \langle I^0_\rho ,A \rangle$ is $\Sigma _m$-elementary by the hypothesis and $\pi \subseteq \tilde \pi$ by (1) -- i.e. $H \cap J^X_\rho \prec _m \langle I^0_\rho , A \rangle$. Since there is a $z \in J^X_\rho$ which satisfies $(\star )$ and $x_i,p \in H \cap J^X_\rho$, there exists such a $z \in H \cap J^X_\rho$.  
\smallskip \\
Let $H \prec _1 \langle I^0_\rho ,A \rangle$. Let $\pi$ be the uncollapse of $H$. Then $\pi$ has a $\Sigma _1$-elementary extension $\tilde \pi : \langle I^0_{\bar \beta},\bar B \rangle \rightarrow \langle I^0_\beta ,B \rangle$. So $H \cong \langle I^0_{\bar \rho},\bar A \rangle$ for some $\bar \rho$ and $\bar A$. I.e. $Cond_A(I^0_\rho )$.
\smallskip \\
(3) $\bar A = \{ \langle i,x \rangle \mid i \in \omega$ $and$ $x \in J^X_{\bar \rho}$ $and$ $\langle I^0_{\bar \beta},\bar B \rangle \models \varphi _i(x, \tilde \pi ^{-1}(p)) \}$
\smallskip \\
Since $\pi : \langle I^0_{\bar \rho}, \bar A \rangle \rightarrow \langle I^0_\rho ,A \rangle$ is $\Sigma _0$-elementary, $\bar A(i,x) \Leftrightarrow A(i,\pi (x))$ for $x \in J^X_{\bar \rho}$. Since $A$ is the standard code of $\langle I^0_\beta , B \rangle$, $A(i,\pi (x))$ $\Leftrightarrow$ $ \langle I^0_\beta , B \rangle \models \varphi _i(\pi (x),p)$. Finally, $\langle I^0_\beta , B \rangle \models \varphi _i(\pi (x),p)$ $ \Leftrightarrow$ $\langle I^0_{\bar \beta},\bar B \rangle \models \varphi _i(x,\tilde \pi ^{-1}(p))$, because $\tilde \pi : \langle I^0_{\bar \beta},\bar B \rangle \rightarrow \langle I^0_\beta , B \rangle$ is $\Sigma _1$-elementary.
\smallskip \\
(4) $\bar \rho$ is the projectum of $\langle I^0_{\bar \beta},\bar B \rangle$
\smallskip \\
By the definition of $H$, $J^X_{\bar \beta}=h_{\bar \beta , \bar B}[\omega \times (J^X_{\bar \rho} \times \{ \tilde \pi ^{-1}(p) \} )]$. So $f(\langle i,x \rangle ) \simeq h_{\bar \beta , \bar B}(i,\langle x, \tilde \pi ^{-1}(p) \rangle )$ is a over $\langle I^0_{\bar \beta},\bar B \rangle$ $\Sigma _1$-definable function such that $f[J^X_{\bar \rho}]=J^X_{\bar \beta}$. It remains to prove that $\langle I^0_{\bar \rho}, C \rangle$ is rudimentary closed for all $C \in \Sigma _1(\langle I^0_{\bar \beta},\bar B \rangle ) \cap \frak{P}(J^X_{\bar \rho })$. By the definition of $H$, there exists an $i\in \omega$ and a $y \in J^X_{\bar \rho}$ such that $x \in C$ $\Leftrightarrow$ $\langle I^0_{\bar \beta},\bar B \rangle \models \varphi _i(\langle x,y \rangle , \tilde \pi ^{-1}(p))$ for all $x \in J^X_{\bar \rho}$. Thus, by (3), $x \in C \Leftrightarrow \bar A(i, \langle x,y \rangle )$. For $u \in J^X_{\bar \rho}$, let $v=\{ \langle i, \langle x,y \rangle \rangle \mid x \in u\}$. Then $v \in J^X_{\bar \rho}$ and $\bar A \cap v \in J^X_{\bar \rho}$, because $\langle I^0_{\bar \rho},\bar A \rangle$ is rudimentary closed by the hypothesis. But $x \in C \cap u$ holds iff $\langle i, \langle x,y \rangle \rangle \in \bar A \cap v$. Finally, $J^X_{\bar \rho}$ is rudimentary closed and therefore  $C\cap u \in J^X_{\bar \rho}$.
\smallskip \\
(5) $\tilde \pi ^{-1}(p)$ is the standard parameter of $\langle I^0_{\bar \beta},\bar B \rangle$
\smallskip \\
By the definition of $H$, $J^X_{\bar \beta}=h_{\bar \beta , \bar B}[\omega \times (J^X_{\bar \rho} \times \{ \tilde \pi ^{-1}(p) \} )]$ and, by (4), $\bar \rho$ is the projectum of $\langle I^0_{\bar \beta},\bar B \rangle$. So we just have to prove that $\tilde \pi ^{-1}(p)$ is the least with this property. Suppose that $\bar p ^\prime < \tilde \pi ^{-1}(p)$ had this property as well. Then there were an $i \in \omega$ and an $x \in J^X_{\bar \rho}$ such that $\tilde \pi ^{-1}(p)=h_{\bar \beta , \bar B}(i,\langle x,\bar p ^\prime \rangle )$.  
Since $\tilde \pi : \langle I^0_{\bar \beta},\bar B \rangle \rightarrow \langle I^0_\beta , B \rangle$ is $\Sigma _1$-elementary, we had $p=h_{\beta ,B}(i,\langle x, p ^\prime \rangle )$ for $p^\prime = \pi ( \bar p ^\prime )<p$. And so also $h_{\beta ,B}[\omega \times (J^X_\rho \times \{ p^\prime \})]=J^X_\beta$. That contradicts the definition of $p$.
\smallskip \\
(6) Uniqueness
\smallskip \\
Assume $\langle I^0_{\bar \beta _0}, \bar B _0 \rangle$ and $\langle I^0_{\bar \beta _1} , \bar B _1 \rangle$ both have $\bar \rho$ as projectum and $\bar A$ as standard code. Let $\bar p _i$  be the standard parameter of $\langle I^0_{\bar \beta _i}, \bar B _i \rangle$. Then, for all $j \in \omega$ and $x \in J^X_{\bar \rho}$, $\langle I^0_{\bar \beta _0}, \bar B _0 \rangle \models \varphi _j( x, \bar p _0)$ iff $\bar A(j,x)$ iff $\langle I^0_{\bar \beta _1} , \bar B _1 \rangle \models \varphi _j( x, \bar p _1)$. So $\sigma (h_{\bar \beta _0,\bar B _0}(j,\langle x, \bar p_0 \rangle )) \simeq h_{\bar \beta _1,\bar B _1}(j,\langle x, \bar p_1 \rangle )$ defines an isomorphism $\sigma : \langle I^0_{\bar \beta _0}, \bar B _0 \rangle \cong \langle I^0_{\bar \beta _0}, \bar B _0 \rangle$, because, for both $i$, $h_{\bar \beta _i,\bar B _i}[\omega \times (J^X_{\bar \rho} \times \{ \bar p _i \} )]=J^X_{\bar \beta _i}$ holds. But since both structures are transitive, $\sigma$ must be the identity. Finally, let $\tilde \pi _0 : \langle I^0_{\bar \beta}, \bar B \rangle \rightarrow \langle I^0_\beta , B \rangle$ and $\tilde \pi _1 : \langle I^0_{\bar \beta}, \bar B \rangle \rightarrow \langle I^0_\beta , B \rangle$ be $\Sigma _1$-elementary extensions of $\pi$. Let $\bar p$ be the standard parameter of $\langle I^0_{\bar \beta} , \bar B \rangle$. Then, for every $y \in J^X_{\bar \beta}$, there is an $x \in J^X_{\bar \rho}$ and a $j \in \omega$ such that $y=h_{\bar \beta , \bar B}(j, \langle x, \bar p \rangle)$ -- and $\tilde \pi _o(y) = h_{\beta , B}(j,\pi (x), \pi (p))=\tilde \pi _1(y)$. Thus $\tilde \pi _0=\tilde \pi _1$. $\Box$ 
\medskip \\
To code the $\Sigma _n$ information of $I_\beta$ where $\beta \in S^X$ in a structure $\langle I^0_\rho , A \rangle$, one iterates this process.
\smallskip \\
For $n \geq 0$, $\beta \in S^X$, let
\smallskip

$\rho ^0=\beta$, $p^0=\emptyset$, $A^0=X_\beta$
\smallskip 

$\rho ^{n+1}=$ the projectum of $\langle I^0_{\rho ^n},A^n \rangle$
\smallskip

$p^{n+1}=$ the standard parameter of $\langle I^0_{\rho ^n},A^n \rangle$
\smallskip

$A^{n+1}=$ the standard code of $\langle I^0_{\rho ^n},A^n \rangle$.
\smallskip \\
Call
\smallskip

$\rho ^n$ the $n$-th projectum of $\beta$,
\smallskip

$p^n$ the $n$-th (standard) parameter of $\beta$,
\smallskip

$A^n$ the $n$-th (standard) code of $\beta$.
\smallskip \\
By lemma 17, $A^n \subseteq J^X_{\rho ^n}$ is $\Sigma _n$-definable over $I_\beta$ and, for all $m \geq 1$,
$$\Sigma _{n+m}(I_\beta ) \cap \frak{P}(J^X_{\rho ^n})=\Sigma _m(\langle I^0_{\rho ^n},A^n \rangle ).$$
From lemma 19, we get by induction:
\smallskip \\
For $\beta > \omega$, $n \geq 1$, $m \geq 0$, let $\rho ^n$ be the $n$-th projectum and $A^n$ be the $n$-th code of $\beta$. Let $\langle I^0_{\bar \rho},\bar A \rangle$ be a rudimentary closed structure and $\pi : \langle I^0_{\bar \rho}, \bar A \rangle \rightarrow \langle I^0_{\rho ^n},A^n \rangle$ be $\Sigma _m$-elementary. Then:
\smallskip \\
(1) There is a unique $\bar \beta \geq \bar \rho$ such that $\bar \rho$ is the $n$-th projectum and $\bar A$ is the $n$-th code of $\bar \beta$.
\smallskip \\
For $k \leq n$ let

$\rho ^k$ be the $k$-th projectum of $\beta$,

$p^k$ the $k$-th parameter of $\beta$,

$A^k$ the $k$-th code of $\beta$
\\
and

$\bar \rho ^k$ the $k$-th projectum of $\bar \beta$,

$\bar p ^k$ the $k$-th parameter of $\bar \beta$,

$\bar A^k$ the $k$-th code of $\bar \beta$.
\smallskip \\
(2) There exists a unique extension $\tilde \pi$ of $\pi$ such that, for all $0 \leq k \leq n$,
\smallskip 

$\tilde \pi \upharpoonright J^X_{\bar \rho ^k} : \langle I^0_{\bar \rho ^k}, \bar A^k \rangle \rightarrow \langle I^0_{\rho ^k}, A^k \rangle$ is $\Sigma _{m+n-k}$-elementary
\smallskip

and $\tilde \pi (\bar p^k)=p^k$.
\medskip \\
{\bf Lemma 20}
\smallskip \\
Let $\omega <\beta \in S^X$. Then all projecta of $\beta$ exist.
\smallskip \\
{\bf Proof} by induction on $n$. That $\rho ^0$ exists is clear. So suppose that the first projecta $\rho ^0, \dots ,\rho ^{n-1}, \rho :=\rho ^n$, the parameters $p^0, \dots , p^n$ and the codes $A^0, \dots A^{n-1},A:=A^n$ of $\beta$ exist. Let $\gamma \in Lim$ be minimal such that there is some over $\langle I^0_\rho ,A \rangle$ $\Sigma _1$-definable function $f$ such that $f[J_\gamma ^X]=J^X_\rho$. Let $C \in \Sigma _1(\langle I^0_\rho ,A \rangle ) \cap \frak{P}(J^X_\gamma )$. We have to prove that $\langle I^0_\gamma ,C \rangle$ is rudimentary closed. If $\gamma = \omega$, then $J^X_\gamma =H_\omega$, and this is obvious. If $\gamma > \omega$, then $\gamma \in Lim^2$ by the definition of $\gamma$. Then it suffices to show $C\cap J^X_\delta \in J^X_\gamma$ for $\delta \in Lim \cap \gamma$. Let $B:=C \cap J^X_\delta$ be definable over $\langle I^0_\rho ,A \rangle$ with parameter $q$. Since obviously $\gamma \leq \rho$, $C \cap J^X_\delta$ is  $\Sigma _n$-definable over $I_\beta$ with parameters $p_1, \dots ,p^n,q$ by lemma 17. So let $\varphi$ be a $\Sigma _n$ formula such that $x \in C$ $\Leftrightarrow$ $I_\beta \models \varphi (x,p^1, \dots , p^n,q)$. Let 

$H_{n+1}:=h_{\rho ^n,A^n}[\omega \times (J^X_\delta \times \{ q\} )]$

$H_n:=h_{\rho ^{n-1},A^{n-1}}[\omega \times (H_n \times \{ p^n\} )]$

$H_{n-1}:=h_{\rho ^{n-2},A^{n-2}}[\omega \times (H_{n-1} \times \{ p^{n-1}\} )]$

etc.
\\
Since $L[X]$ has condensation, there is an $I_\mu$ such that $H_1 \cong I_\mu$. Let $\pi$ be the uncollapse of $H_1$. Then $\pi$ is the extension of the collapse of $H_{n+1}$ defined in the proof of lemma 19. Therefore it is $\Sigma _{n+1}$-elementary. Since $B \subseteq J^X_\delta$ and $\pi \upharpoonright J^X_\delta =id \upharpoonright J^X_\delta$, we get $x \in B$ $\Leftrightarrow$ $I_\mu \models \varphi (x, \pi ^{-1}(p^1), \dots , \pi ^{-1}(p^n), \pi ^{-1}(q))$. So $B$ is indeed already $\Sigma _n$-definable over $I_\mu$. Thus $B \in J^X_{\mu +1}$ by lemma 8. But now we are done because $\mu < \rho$. For, if

$h_{n+1}(\langle i,x \rangle )=h_{\rho ^n,A^n}(i,\langle x,p \rangle )$

$h_n(\langle i,x \rangle )=h_{\rho ^{n-1},A^{n-1}}(i,\langle x,p^n\rangle )$

etc.
\\
then the function $h=h_1 \circ \dots \circ h_{n+1}$ is $\Sigma _{n+1}$-definable over $I_\beta$. Thus the function $\bar h =\pi [h \cap (H_1 \times H_1)]$ is $\Sigma _{n+1}$-definable over $I_\mu$ and $\bar h [J^X_\delta ]=J^X_\mu$. So  $\bar h \cap (J^X_\rho )^2$ is $\Sigma _1$-definable over $\langle I^0_\rho ,A \rangle$ by lemma 17 and lemma 19. And by the definition of $\gamma$, there is an over $\langle I^0_\rho , A \rangle$ $\Sigma _1$-definable function $f$ such that $f[J^X_\gamma ]=J^X_\rho$. So if we had $\mu \geq \rho$, then $f \circ \bar h$ was an over $\langle I^0_\rho , A\rangle$ $\Sigma _1$-definable function such that $(f \circ \bar h)[J^X_\delta ]=J^X_\rho$. That contradicts the minimality of $\gamma$. $\Box$
\medskip \\
Let $\omega <\nu \in S^X$, $\rho ^n$ the $n$-th projectum of $\nu$, $p^n$ the $n$-th parameter and $A^n$ the $n$-th Code. Let

$h_{n+1}(i,x)=h_{\rho ^n,A^n}(i,x)$

$h_n(\langle i,x \rangle )=h_{\rho ^{n-1},A^{n-1}}(i,\langle x,p^n\rangle )$

etc.
\\
Then define $$h_\nu ^{n+1}=h_1 \circ \dots \circ h_{n+1}.$$
We have:
\smallskip\\
(1) $h_\nu ^n$ is $\Sigma _n$-definable over $I_\nu$
\smallskip\\
(2) $h^n_\nu [\omega \times Q] \prec_n I_\nu$, if $Q \subseteq J^X_{\rho ^{n-1}}$ is closed under ordered pairs. 
\medskip \\
{\bf Lemma 21}
\smallskip \\
Let $\omega < \beta \in S^X$ and $n \geq 1$. Then
\smallskip \\
(1) the least ordinal $\gamma \in Lim$ such that there is a over $I_\beta$ $\Sigma _n$-definable function $f$ such that $f[J^X_\gamma]=J^X_\beta$,
\smallskip \\
(2) the last ordinal $\gamma \in Lim$ such that $\langle I^0_\gamma,C \rangle$ is rudimentary closed for all $C \in \Sigma _n(I_\beta ) \cap \frak{P}(J^X_\gamma )$,
\smallskip \\
(3) the least ordinal $\gamma \in Lim$ such that $\frak{P}(\gamma ) \cap \Sigma _n(I_\beta ) \nsubseteq J^X_\beta$,
\smallskip \\
is the $n$-th projectum of $\beta$.
\smallskip \\
{\bf Proof:}
\smallskip \\
(1) By the definition of the $n$-th projectum, there is an over $\langle I^0_{\rho ^{n-1}},A^{n-1}\rangle$ $\Sigma _1$-definable $f^n$ such that $f^n[J^X_{\rho ^n}]=J^X_{\rho ^{n-1}}$, an over  $\langle I^0_{\rho ^{n-2}},A^{n-2}\rangle$ $\Sigma _1$-definable $f^{n-1}$ such that $f^{n-1}[J^X_{\rho ^{n-1}}]=J^X_{\rho ^{n-2}}$, etc. But then $f^k$ is $\Sigma _k$-definable over $I_\beta$ by lemma 17. So $f=f^1 \circ f^2 \circ \dots \circ f^n$ is $\Sigma _n$-definable over $I_\beta$ and $f[J^X_{\rho ^n}]=J^X_\beta$.

On the other hand, the projectum $\bar \rho$ of a rudimentary closed structure $\langle I^0_\beta , B \rangle$ is the least $\bar \rho$ such that there is an over $\langle I^0_\beta ,B \rangle$ $\Sigma _1$-definable function $f$ such that $f[J^X_{\bar \rho}]=J^X_\beta$. For, suppose there is no such $\rho < \bar \rho$ such that such an  $f$, $f[J^X_\rho ]=J^X_\beta$, exists. Then the proof of lemma 16 provides a contradiction. So if there was a $\gamma < \rho ^n$ such that there is an over $I_\beta$ $\Sigma _n$-definable function $f$ such that $f[J^X_\gamma ]=J^X_\beta$, then $g:=f \cap (J^X_{\rho ^{n-1}})^2$ would be an over $\langle I^0_{\rho ^{n-1}},A^{n-1}\rangle$ $\Sigma _1$-definable function such that $g[J^X_\gamma ]=J^X_{\rho ^{n-1}}$. But this is impossible.
\smallskip \\
(2) By the definition of the $n$-th projectum, $\langle I^0_{\rho ^n}, C \rangle$ is rudimentary closed for all $C \in \Sigma _1(\langle I^0_{\rho ^{n-1}}, A ^{n-1} \rangle ) \cap \frak{P}(J^X_{\rho ^n})$. But by lemma 17, $\Sigma _1(\langle I^0_{\rho ^{n-1}},A^{n-1} \rangle )= \Sigma _n(I_\beta ) \cap \frak{P}(J^X_{\rho ^{n-1}})$. So, since $\rho ^n \leq \rho ^{n-1}$, $\langle I^0_{\rho ^n}, C \rangle$ is rudimentary closed for all $C \in \Sigma _n(I_\beta ) \cap \frak{P}(J^X_{\rho ^n})$.

Assume $\gamma$ were a larger ordinal $ \in Lim$ having this property. Let $f$ be, by (1), an over $I_\beta$ $\Sigma _n$-definable function such that $f[J^X_{\rho ^n}]=J^X_\beta$. Set $C=\{u \in J^X_{\rho ^n} \mid u \notin f(u) \}$. Then $C$ is $\Sigma _n$-definable over $I_\beta$ and $C \subseteq J^X_{\rho ^n}$. So $\langle J^X_\gamma , C \rangle$ was rudimentary closed. And therefore $C=C \cap J^X_{\rho ^n} \in J^X_\gamma \subseteq J^X_\beta$ and $C=f(u)$ for some $u \in J^X_{\rho ^n}$. But this implies the contradiction that $u \in f(u) \Leftrightarrow u \in C \Leftrightarrow u \notin f(u)$.
\smallskip \\
(3) Let $\rho := \rho ^n$ and $f$ by (1) an over $I_\beta$ $\Sigma _n$-definable function such that $f[J^X_\rho ]=J^X_\beta$. Let $j$ be an over $I^0_\rho$ $\Sigma _1$-definable function from $\rho$ onto $J^X_\rho$. Let $C=\{ \nu \in \rho \mid \nu \notin f \circ j(\nu )\}$. Then $C$ is an over $I_\beta$ $\Sigma _n$-definable subset of  $\rho$. If $C \in J^X_\beta$, then there would be a $\nu \in \rho$ such that $C=f\circ j(\nu )$, and we had the contradiction $\nu \in C \Leftrightarrow \nu \notin f \circ j(\nu ) \Leftrightarrow \nu \notin C$. Thus $\frak{P}(\rho ) \cap \Sigma _n(I_\beta ) \nsubseteq J^X_\beta$. But if $\gamma \in Lim \cap \rho$ and $D \in \frak{P}(\gamma ) \cap \Sigma _n(I_\beta )$, then $D=D \cap J^X_\gamma \in J^X_\rho \subseteq J^X_\beta$. So $\frak{P}(\gamma ) \cap \Sigma _n(I_\beta ) \subseteq J^X_\beta$. $\Box$

\section{Morasses}
Let $\omega _1 \leq \beta$, $S=Lim \cap \omega _{1+\beta}$ and $\kappa :=\omega _{1+\beta}$. 
\smallskip\\
 We write $Card$ for the class of cardinals and $RCard$ for the class of regular cardinals.
\smallskip\\
Let $\vartriangleleft$ be a binary relation on $S$ such that:
\smallskip \\
(a) If $\nu \vartriangleleft \tau$, then $\nu < \tau$.

For all $\nu \in S - RCard$, $\{ \tau \mid \nu \vartriangleleft \tau \}$ is closed. 

For $\nu \in S - RCard$, there is a largest $\mu$ such that $\nu \trianglelefteq \mu$.
\smallskip \\
Let $\mu _\nu$ be this largest $\mu$ with $\nu \trianglelefteq \mu$.
\smallskip\\
Let 
$$\nu \sqsubseteq \tau : \Leftrightarrow \nu \in Lim(\{ \delta \mid \delta \vartriangleleft \tau \} ) \cup \{ \delta \mid \delta \trianglelefteq \tau \} .$$
(b) $\sqsubseteq$ is a (many-rooted) tree. 
\smallskip\\
Hence, if $\nu \notin RCard$ is a successor in $\sqsubset$, then $\mu _\nu$ is the largest $\mu$ such that $\nu \sqsubseteq \mu$. To see this, let $\mu ^\ast _\nu$ be the largest $\mu$ such that $\nu \sqsubseteq \mu$. It is clear that $\mu _\nu \leq \mu _\nu ^\ast$, since $\nu \trianglelefteq \mu$ implies $\nu \sqsubseteq \mu$. So assume that $\mu _\nu < \mu _\nu ^\ast$. Then $\nu \not\vartriangleleft \mu_\nu^\ast$ by the definition of $\mu _\nu$. Hence $\nu \in Lim(\{ \delta \mid \delta \vartriangleleft \mu _\nu ^\ast\} )$ and $\nu \in Lim(\{ \delta \mid \delta \sqsubseteq \mu _\nu^\ast\} )$. Therefore, $\nu \in Lim(\sqsubseteq)$ since $\sqsubseteq$ is a tree. That contradicts our assumption that $\nu$ is a successor in $\sqsubset$.
\smallskip \\
For $\alpha \in S$, let $|\alpha|$ be the rank of $\alpha$ in this tree. Let 
\smallskip 

$S^+ := \{ \nu \in S \mid \nu$ is a successor in $\sqsubset \}$
\smallskip

$S^0 := \{ \alpha \in S \mid |\alpha| =0\}$
\smallskip

$\widehat{S^+}:=\{ \mu _\tau \mid \tau \in S^+ - RCard \}$
\smallskip

$\widehat{S}:=\{ \mu _\tau \mid \tau \in S -RCard \}$.
\smallskip \\
Let $S_\alpha := \{ \nu \in S \mid \nu$ is a direct successor of $\alpha$ in $\sqsubset \}$. For $\nu \in S^+$, let $\alpha _\nu$ be the direct predecessor of $\nu$ in $\sqsubset$. For $\nu \in S^0$, let $\alpha _\nu :=0$. For $\nu \not\in S^+ \cup S^0$, let $\alpha _\nu :=\nu$.
\smallskip \\
(c) For $\nu ,\tau \in (S^+ \cup S^0)-RCard$ such that $\alpha _\nu =\alpha _\tau$, suppose:
$$\nu < \tau \quad \Rightarrow \quad \mu _\nu < \tau .$$
For all $\alpha \in S$, suppose:
\begin{tabbing}
(d) \= $S_\alpha$ is closed \\[0.25ex]
(e) \> $card(S_\alpha ) \leq \alpha ^+$ \\[0.25ex]
\> $card(S _\alpha ) \leq card(\alpha )$ if $card(\alpha ) < \alpha$ \\[0.25ex]
(f) \> $\omega _1 = max(S^0 ) = sup(S^0 \cap \omega _1 )$ \\[0.25ex]
\> $\omega _{1+i+1}=max(S_{\omega _{1+i}})=sup(S_{\omega _{1+i}} \cap \omega _{1+i+1})$ for all $i<\beta$. \\
\end{tabbing}
Let $D=\langle D_\nu \mid \nu \in \widehat{S} \rangle$ be a sequence such that $D_\nu \subseteq J^D_\nu$. To simplify matters, my definition of $J^D_\nu$ is such that $J^D_\nu \cap On =\nu$ (see section 3 or [SchZe]).
\medskip \\
Let an $\langle S,\vartriangleleft ,D\rangle$-maplet $f$ be a triple $\langle \bar \nu ,|f| , \nu \rangle$ such that $\bar \nu ,\nu \in S - RCard$ and $|f| :J^D_{\mu _{\bar \nu}} \rightarrow J^D_{\mu _\nu}$.
\smallskip \\
Let $f=\langle \bar \nu ,|f| , \nu \rangle$ be an  $\langle S,\vartriangleleft ,D\rangle$-maplet. Then we define $d(f)$ and $r(f)$ by $d(f)=\bar \nu$ and $r(f)=\nu$. Set $f(x):=|f| (x)$ for $x \in J^D_{\mu _{\bar \nu}}$ and $f(\mu_{\bar\nu} ):=\mu _\nu$. But $dom(f)$, $rng(f)$, $f \upharpoonright X$, etc. keep their usual set-theoretical meaning, i.e. $dom(f)=dom(|f| )$, $rng(f)=rng(|f| )$, $f \upharpoonright X = |f| \upharpoonright X$, etc. 
\smallskip \\
For $\bar \tau \leq \mu_{ \bar \nu}$, let $f^{( \bar \tau )}=\langle \bar \tau , |f| \upharpoonright J^D_{\mu _{\bar \tau}}, \tau \rangle$ where $\tau =f( \bar \tau )$. Of course, $f^{( \bar \tau )}$ needs not to be a maplet. The same is true for the following definitions.  Let $f^{-1}=\langle \nu , |f| ^{-1}, \bar \nu \rangle$. For $g= \langle \nu , |g| , \nu ^\prime \rangle$ and $f=\langle \bar \nu ,|f| , \nu \rangle$, let $g \circ f= \langle \bar \nu , | g| \circ |f| , \nu ^\prime \rangle$. If $g= \langle \nu ^\prime , | g | , \nu \rangle$ and $f=\langle \bar \nu ,| f| , \nu \rangle$ such that $rng(f) \subseteq rng(g)$, then set $g^{-1}f=\langle \bar \nu , | g| ^{-1} \mid f \mid , \nu ^\prime \rangle$. Finally set $id _\nu = \langle \nu , id \upharpoonright J^D_{\mu _\nu} , \nu \rangle$.
\smallskip \\
Let $\frak{F}$ be a set of  $\langle S,\vartriangleleft ,D\rangle$-maplets $f=\langle \bar \nu ,|f| , \nu \rangle$ such that the following \\
holds:
\smallskip \\
(0) $f(\bar \nu )=\nu$, $f(\alpha _{\bar \nu})=\alpha _\nu$ and $|f|$ is order-preserving. 
\smallskip \\
(1) For $f \neq id _{\bar \nu}$, there is some $\beta \sqsubseteq \alpha _{\bar \nu}$ such that $f \upharpoonright \beta = id \upharpoonright \beta$ and $f(\beta )>\beta$.
\smallskip \\
(2) If $\bar \tau \in S^+$ and $\bar \nu \sqsubset \bar \tau \sqsubseteq \mu_{\bar \nu}$, then $f^{(\bar \tau )}\in \frak{F}$.
\smallskip \\
(3) If $f,g \in \frak{F}$ and $d(g)=r(f)$, then $g \circ f \in \frak{F}$.
\smallskip\\
(4) If $f,g \in \frak{F}$, $r(g)=r(f)$ and $rng(f)\subseteq rng(g)$, then $g^{-1} \circ f \in \frak{F}$.
\smallskip\\
We write $f: \bar \nu \Rightarrow \nu$ if $f=\langle \bar \nu , |f| , \nu \rangle \in \frak{F}$. If $f \in \frak{F}$ and $r(f)=\nu$, then we write $f \Rightarrow \nu$. The uniquely determined $\beta$ in (1) shall be denoted by $\beta (f)$.
\smallskip \\
Say $f \in \frak{F}$ is minimal for a property $P(f)$ if $P(f)$ holds and $P(g)$ implies $g^{-1}f \in \frak{F}$. 
\smallskip \\
Let 
\smallskip

$f_{(u,x ,\nu )}=$ the unique minimal $f \in \frak{F}$ for $f \Rightarrow \nu$ and $u \cup \{ x \} \subseteq rng(f)$,
\smallskip \\
if such an $f$ exists. The axioms of the morass will guarantee that $f_{(u,x ,\nu )}$ always exists if $\nu \in S-RCard^{L_\kappa[D]}$. Therefore, we will always assume and explicitly mention that $\nu \in S-RCard^{L_\kappa[D]}$ when $f_{(u,x ,\nu )}$ is mentioned. 
\smallskip \\
Say $\nu \in S-RCard^{L_\kappa[D]}$ is independent if $d(f_{(\beta ,0,\nu )}) < \alpha _\nu$ holds for all $\beta < \alpha _\nu$.
\smallskip \\
For $\tau \sqsubseteq \nu\in S-RCard^{L_\kappa[D]}$, say $\nu$ is $\xi$-dependent on $\tau$ if $f_{(\alpha _\tau , \xi , \nu )} = id _\nu$.
\medskip \\
For $f \in \frak{F}$, let $\lambda (f):= sup(f[d(f)])$.
\medskip\\
For $\nu\in S-RCard^{L_\kappa[D]}$ let
$$C_\nu =\{ \lambda (f) <\nu \mid f \Rightarrow \nu \} $$
$$\Lambda (x ,\nu )=\{ \lambda (f_{(\beta , x , \nu )}) < \nu \mid \beta < \nu \} .$$
It will be shown that $C_\nu$ and $\Lambda (x ,\nu )$ are closed in $\nu$.
\smallskip\\
Recursively define a function $q_\nu:k_\nu +1 \rightarrow On$, where $k_\nu \in \omega$:
\smallskip 

$q_\nu (0)=0$
\smallskip

$q_\nu (k+1) = max(\Lambda (q_\nu \upharpoonright (k+1) ,\nu ))$
\smallskip \\
if $max(\Lambda (q_\nu \upharpoonright (k+1) ,\nu ))$ exists. The axioms will guarantee that this recursion breaks off (see lemma 4 below), i.e. there is some $k_\nu$ such that either
\smallskip

$\Lambda (q_\nu \upharpoonright ({k_\nu}+1),\nu ) =\emptyset$
\smallskip \\
or
\smallskip

$\Lambda (q_\nu \upharpoonright ({k_\nu}+1),\nu )$ is unbounded in $\nu$.
\medskip \\
Define by recursion on $1 \leq n \in \omega$, simultaneouly for all $\nu \in S-RCard^{L_\kappa[D]}$,  $\beta \in \nu$ and $x \in J^D_{\mu_\nu}$ the following notions:
\smallskip

$f^1_{(\beta ,x,\nu )}=f_{(\beta ,x ,\nu )}$
\smallskip

$\tau (n,\nu )$ $=$ the least $\tau \in S^0 \cup S^+ \cup \widehat{S}$ such that for some $x \in J^D_{\mu _\nu}$
$$f^n_{(\alpha _\tau ,x ,\nu )}=id_\nu$$

$x (n,\nu )$ $=$ the least $x \in J^D_{\mu _\nu}$ such that $f^n_{(\alpha _{\tau (n,\nu )} ,x ,\nu )}=id_\nu$
\smallskip 

$K^n_\nu =\{ d(f^n_{(\beta ,x (n,\nu ),\nu )})< \alpha _{\tau (n,\nu )} \mid \beta < \nu \}$
\smallskip

$f \Rightarrow _n \nu$ iff $f \Rightarrow \nu$ and for all $1 \leq m < n$ $$rng(f) \cap J^D_{\alpha _{\tau (m,\nu )}} \prec _1 \langle J^D_{\alpha _{\tau (m,\nu )}},D \upharpoonright {\alpha _{\tau (m,\nu )}},K^m_\nu \rangle $$
$$x (m,\nu ) \in rng(f)$$

$f^n_{(u ,\nu )}$ $=$ the minimal $f \Rightarrow _n \nu$ such that $u \subseteq rng(f)$
\smallskip 

$f^n_{(\beta , x ,\nu )}$ $=$ $f^n_{(\beta \cup \{ x \} ,\nu )}$
\smallskip 

$f:\bar \nu \Rightarrow _n \nu$ $:\Leftrightarrow$ $f \Rightarrow _n \nu$ $and$ $f:\bar \nu \Rightarrow \nu$.
\medskip \\
Here definitions are to be understood in Kleene's sense, i.e., that the left side is defined iff the right side is, and in that case, both are equal.
\medskip\\
Let
\smallskip

$n_\nu$ $=$ the least $n$ such that $f^n_{(\gamma ,x ,\mu _\nu )}$ is confinal in $\nu$ for some $x \in J^D_{\mu _\nu}$, $\gamma \sqsubset \nu$
\smallskip

$x _\nu$ $=$ the least $x$ such that $f^{n_\nu}_{(\alpha _\nu ,x ,\mu _\nu )}=id_{\mu _\nu}$.
\smallskip \\
Let
\smallskip 

$\alpha ^\ast _\nu =\alpha _\nu$ if $\nu \in S^+$
\smallskip 

$\alpha _\nu ^\ast = sup \{ \alpha < \nu \mid \beta (f^{n_\nu}_{(\alpha ,x _\nu ,\mu _\nu )})=\alpha \}$ if $\nu \notin S^+$.
\smallskip\\
Let $P_\nu := \{ x _\tau \mid \nu \sqsubset \tau \sqsubseteq \mu _\nu ,\tau \in S^+ \} \cup \{ x _\nu \}$.
\bigskip \\
We say that $\frak{M} = \langle S,\vartriangleleft ,\frak{F},D \rangle$ is an $(\omega _1 , \beta )$-morass if the following axioms hold:
\medskip \\
{\bf (MP -- minimum principle)}
\smallskip\\
If $\nu \in S-RCard^{L_\kappa [D]}$ and $x \in J^D_{\mu _\nu}$, then $f_{(0,x ,\nu)}$ exists.
\medskip \\
{\bf (LP1 -- first logical preservation axiom)}
\smallskip\\
If $f:\bar \nu \Rightarrow \nu$, then $|f| : \langle J^D_{\mu _{\bar \nu}} , D \upharpoonright {\mu _{\bar \nu}} \rangle \rightarrow \langle J^D_{\mu _\nu} , D \upharpoonright {\mu _\nu} \rangle$ is $\Sigma _1$-elementary.
\medskip\\
{\bf (LP2 -- second logical preservation axiom)}
\smallskip\\
Let $f:\bar \nu \Rightarrow \nu$ and $f(\bar x )=x$. Then $$(f \upharpoonright J^D_{\bar \nu }):\langle J^D_{\bar \nu} ,D \upharpoonright {\bar \nu},\Lambda (\bar x , \bar \nu ) \rangle \rightarrow \langle J^D_\nu ,D \upharpoonright \nu ,\Lambda (x ,\nu ) \rangle$$ is $\Sigma _0$-elementary.
\medskip \\
{\bf (CP1 -- first continuity principle)}
\smallskip\\
For $i \leq j < \lambda$, let $f_i:\nu _i \Rightarrow \nu$ and $g_{ij}:\nu _i \Rightarrow \nu _j$ such that $g_{ij}=f^{-1}_jf_i$. Let $\langle g_i \mid i<\lambda \rangle$ be the transitive, direct limit of the directed system $\langle g_{ij} \mid i \leq j<\lambda \rangle$ and $hg_i= f_i$ for all $i<\lambda$. Then $g_i,h \in \frak{F}$.
\pagebreak\\
{\bf (CP2 -- second continuity principle)}
\smallskip\\
Let $f:\bar \nu \Rightarrow \nu$ and $\lambda =sup(f[\bar \nu ])$. If, for some $\bar \lambda$, $h:\langle J^{\bar D}_{\bar \lambda} , \bar D \rangle \rightarrow \langle J^D_\lambda, D \upharpoonright \lambda \rangle$ is $\Sigma _1$-elementary and $rng(f \upharpoonright J^D_{\bar \nu} ) \subseteq rng(h)$, then there is some $g:\bar \lambda \Rightarrow \lambda$ such that $g \upharpoonright J^{\bar D}_{\bar \lambda} = h$. 
\medskip \\
{\bf (CP3 -- third continuity principle)}
\smallskip \\
If $C_\nu =\{ \lambda (f) < \nu \mid f \Rightarrow \nu \}$ is unbounded in $\nu\in S-RCard^{L_\kappa[D]}$, then the following holds for all $x \in J^D_{\mu _\nu}$:
$$rng(f_{(0,x ,\nu )})=\bigcup \{ rng(f_{(0,x ,\lambda )}) \mid \lambda \in C_\nu \} .$$
{\bf (DP1 -- first dependency axiom)}
\smallskip\\
If $\mu _\nu < \mu _{\alpha _\nu}$, then $\nu\in S-RCard^{L_\kappa[D]}$ is independent.
\medskip \\
{\bf (DP2 -- second dependency axiom)}
\smallskip \\
If $\nu\in S-RCard^{L_\kappa[D]}$ is $\eta$-dependent on $\tau \sqsubseteq \nu$, $\tau \in S^+$, $f:\bar \nu \Rightarrow \nu$, $f(\bar \tau )=\tau$ and $\eta \in rng(f)$, then $f^{(\bar \tau )}: \bar \tau \Rightarrow \tau$.
\medskip \\
{\bf (DP3 -- third dependency axiom)}
\smallskip\\
For $\nu \in \widehat{S}-RCard^{L_\kappa[D]}$ and $1 \leq n \in \omega$, the following holds:
\smallskip \\
(a) If $f^n_{(\alpha _\tau , x ,\nu )}=id_\nu$, $\tau \in S^+ \cup S^0$ and $\tau \sqsubseteq \nu$, then $\mu _\nu =\mu _\tau$.
\smallskip \\
(b) If $\beta < \alpha _{\tau (n,\nu )}$, then also $d(f^n_{(\beta ,x (n,\nu ),\nu )})< \alpha _{\tau (n,\nu )}$.
\medskip \\
{\bf (DF -- definability axiom)}
\smallskip \\
(a) If $f_{(0,z_0 ,\nu )}=id_\nu$ for some $\nu \in \widehat{S}-RCard^{L_\kappa[D]}$ and $z_0 \in J^D_{\mu _\nu}$, then
$$\{ \langle z ,x,f_{(0,z ,\nu )}(x) \rangle \mid z \in J^D_{\mu _\nu}, x \in dom(f_{(0,z ,\nu )}) \}$$
is uniformly  definable over $\langle J^D_{\mu _\nu},D\upharpoonright {\mu _\nu},D_{\mu _\nu} \rangle$.
\smallskip \\
(b) For all $\nu \in S-RCard^{L_\kappa [D]}$, $x \in J^D_{\mu _\nu}$, the following holds:
$$f_{(0,x ,\nu )}=f^{n_\nu}_{(0,\langle x ,\nu , \alpha ^\ast _\nu ,P_\nu \rangle ,\mu _\nu )}.$$
\smallskip \\
This finishes the definition of an $(\omega_1,\beta)$-morass.
\bigskip\\
A consequence of the axioms is ($\times$) by [Irr2]::
\bigskip\\
{\bf Theorem}
\smallskip\\
$$\{ \langle z ,\tau ,x,f_{(0,z ,\tau )}(x) \rangle \mid  \tau < \nu ,\mu _\tau =\nu , z \in J^D_{\mu _\tau}, x \in dom(f_{(0,z ,\tau )}) \}$$
$$\cup \{ \langle z ,x,f_{(0,z ,\nu )}(x) \rangle \mid \mu _\nu =\nu ,z \in J^D_{\mu _\nu}, x \in dom(f_{(0,z ,\nu )}) \}$$
$$\cup (\sqsubset \cap \nu ^2)$$
is for all $\nu \in S$ uniformly definable over $\langle J^D_\nu ,D\upharpoonright \nu ,D_\nu \rangle$.
\bigskip\\
A structure  $\frak{M} = \langle S,\vartriangleleft ,\frak{F},D \rangle$ is called an $\omega _{1+\beta}$-standard morass if it satifies all axioms of an $( \omega _1 , \beta )$-morass except {\bf (DF)} which is replaced by:
\medskip 

$\nu \vartriangleleft \tau$ $\Rightarrow$ $\nu$ is regular in $J^D_\tau$
\medskip \\
and there are functions $\sigma _{(x ,\nu )}$ for $\nu \in \widehat{S}$ and $x \in J^D_\nu$ such that:
\medskip \\
{\bf (MP)$^+$}
\smallskip \\
$\sigma _{(x ,\nu )}[\omega ]=rng(f_{(0,x ,\nu )})$
\medskip \\
{\bf (CP1)$^+$}
\smallskip \\
If $f:\bar \nu \Rightarrow \nu$ and $f(\bar x )=x$, then $\sigma _{(x ,\nu )}=f \circ \sigma _{(\bar x ,\bar \nu )}$.
\medskip \\
{\bf (CP3)$^+$}
\smallskip \\
If $C_\nu$ is unbounded in $\nu$, then $\sigma _{(x ,\nu )}=\bigcup \{ \sigma _{(x ,\lambda )} \mid \lambda \in C_\nu , x \in J^D_\lambda \}$.
\medskip \\
{\bf (DF)$^+$}
\smallskip \\
(a) If $f_{(0,x ,\nu )}=id_\nu$ for some $x \in J^D_\nu$, then
$$\{ \langle i ,z ,\sigma _{(z ,\nu )}(i) \rangle \mid z \in J^D_\nu, i \in dom(\sigma _{(z ,\nu )}) \}$$
is uniformly definable over $\langle J^D_{\mu _\nu},D\upharpoonright {\mu _\nu},D_{\mu _\nu} \rangle$.
\smallskip \\
(b) If $C_\nu$ is unbounded in $\nu$, then $D_\nu =C_\nu$. If it is bounded, then $D_\nu =\{ \langle i,\sigma _{(q_\nu ,\nu )} (i) \rangle \mid i \in dom (\sigma _{(q_\nu ,\nu )}) \}$. 
\bigskip \\
Now, I am going to construct a $\kappa$-standard morass. 
\smallskip\\
Let $\beta (\nu )$ be the least $\beta$ such that $J^X_{\beta +1}\models \nu$ singular.
\smallskip\\
Let $L_\kappa [X]$ satisfy amenability, condensation and coherence such that $S^X=\{ \beta (\nu ) \mid \nu$ singular in $L_\kappa [X]\}$ and $Card^{L_\kappa [X]}=Card\cap \kappa$. 
\medskip \\
Let 
$$\nu \vartriangleleft \tau \quad : \Leftrightarrow \quad \nu \quad \hbox{regular in} \quad I_\tau .$$
Let 
$$E=Lim-RCard^{L_\kappa [X]} .$$
For $\nu \in E$, let
\smallskip 

$\beta (\nu )=$ the least $\beta$ such that there is a cofinal $f:a \rightarrow \nu \in Def(I_\beta )$ and $a \subseteq \nu ^{\prime} < \nu$ 
\smallskip 

$n(\nu )=$ the least $n \geq 1$ such that such an $f$ is $\Sigma _n$-definable over $I_{\beta (\nu )}$ 
\smallskip 

$\rho (\nu )=$ the $(n(\nu )-1)$-th projectum of $I_{\beta (\nu )}$
\smallskip 

$A_\nu =$ the $(n(\nu )-1)$-th standard code of $I_{\beta (\nu )}$
\smallskip

$\gamma (\nu )=$ the $n(\nu )$-th projectum of $I_{\beta (\nu )}$. 
\medskip \\
If $\nu \in S^+ -Card$, then the $n(\nu )$-th projectum $\gamma$ of $\beta (\nu )$ is less or equal $\alpha _\nu :=$ the largest cardinal in $I_\nu$: Since $\alpha _\nu$ is the largest cardinal in $I_\nu$, there is, by definition of $\beta (\nu )$ and $n(\nu )$, some  over $I_{\beta (\nu )}$ $\Sigma _{n(\nu )}$-definable function $f$ such that $f[\alpha _\nu ]$ is cofinal in $\nu$. But, since $\nu$ is regular in $\beta (\nu )$, $f$ cannot be an element of $J^X_{\beta (\nu )}$. So $\frak{P}(\nu \times \nu ) \cap \Sigma _{n(\nu )}(I_{\beta (\nu )}) \nsubseteq J^X_{\beta (\nu )}$. By lemma 14, also $\frak{P}(\nu ) \cap \Sigma _{n(\nu )}(I_{\beta (\nu )}) \nsubseteq J^X_{\beta (\nu )}$. Using lemma 21  (3), we get $\gamma \leq \nu$. I.e. there is an over $I_{\beta (\nu )}$ $\Sigma _{n(\nu )}$-definable function $g$ such that $g[\nu ]=J^X_{\beta (\nu )}$. On the other hand, there is, for every $\tau < \nu$ in $J^X_\nu$, a surjection from $\alpha _\nu$ onto $\tau$, because $\alpha _\nu$ is the largest cardinal in $I_\nu$. Let $f_\tau$ be the $<_\nu$-least such. Define $j_1(\sigma , \tau )=f_{f(\tau )}(\sigma )$ for $\sigma , \tau < \nu$. Then $j_1$ is $\Sigma _{n(\nu )}$-definable over $I_{\beta (\nu )}$ and $j_1[\alpha _\nu \times \alpha _\nu ]=\nu$. By lemma 15, we obtain an over $I_{\beta (\nu )}$ $\Sigma _{n(\nu )}$-definable function $j_2$ from a subset of $\alpha _\nu$ onto $\nu$. Thus $g \circ j_2$ is an over $I_{\beta (\nu )}$ $\Sigma _{n(\nu )}$-definable map such that $g \circ j_2[\alpha _\nu ]=J^X_{\beta (\nu )}$.
\smallskip \\
Moreover, $\alpha _\nu < \nu \leq \rho (\nu )$: By definition of $\rho (\nu )$, there is an over $I_{\beta (\nu )}$ $\Sigma _{n(\nu )-1}$-definable function $f$ such that $f[\rho (\nu )]=\beta (\nu )$ if $n(\nu )>1$. But $\nu$ is $\Sigma _{n(\nu )-1}$-regular over $I_{\beta (\nu )}$. Thus $\nu \leq \rho (\nu )$. If $n(\nu )=1$, then $\rho (\nu )=\beta (\nu ) \geq \nu$.
\smallskip\\
By the first inequality, there is a $q$ such that every $x \in J^X_{\rho (\nu )}$ is $\Sigma _1$-definable in $\langle I^0_{\rho (\nu )},A_\nu \rangle$ with parameters from $\alpha _\nu \cup \{ q \}$. Let $p_\nu$ be the $<_{\rho (\nu )}$-least such.
\smallskip \\
Obviously, $p_\tau \leq p_\nu$ if $\nu \sqsubseteq \tau \sqsubseteq \mu _\nu$.
\smallskip \\
Thus $P_\nu := \{ p_\tau \mid \nu \sqsubseteq \tau \sqsubseteq \mu _\nu ,\tau \in S^+ \}$ is finite.
\smallskip \\
Now, let $\nu \in E-S^+$. By definition of $\beta (\nu )$, there exists no cofinal $f:a \rightarrow \nu$ in $J^X_\beta$ such that $a \subseteq \nu ^{\prime} < \nu$. So $\frak{P}(\nu \times \nu ) \cap \Sigma _{n(\nu )}(I_{\beta (\nu )}) \not\subseteq J^X_{\beta (\nu )}$. Then, by lemma 14, $\frak{P}(\nu ) \cap \Sigma _{n(\nu )}(I_{\beta (\nu )}) \not\subseteq J^X_{\beta (\nu )}$. Hence, by lemma 21 (3),
$$\gamma (\nu ) \leq \nu .$$
Assume $\rho (\nu )< \nu$. Then there was an over $I_{\beta (\nu )}$ $\Sigma _{n(\nu )-1}$-definable $f$ such that $f[\rho (\nu )]=\nu$. But this contradicts the definition of $n(\nu )$. So 
$$\nu \leq \rho (\nu ).$$
Using lemma 21 (1), it follows from the first inequality that there is some over $I_{\beta (\nu )}$ $\Sigma _{n(\nu )}$-definable function $f$ such that $f[J^X_\nu ]=J^X_{\beta (\nu )}$. So there is a $p \in J^X_{\rho (\nu )}$ such that every $x \in J^X_{\rho (\nu )}$ is $\Sigma _1$-definable in $\langle I^0_{\rho (\nu )},A_\nu \rangle$ with parameters from $\nu \cup \{ p \}$. Let $p_\nu$ be the least such. 
\smallskip \\
Let
$$\alpha ^\ast _\nu =sup\{\alpha < \nu \mid h_{\rho (\nu ),A_\nu}[\omega \times (J^X_\alpha \times \{ p_\nu \} )] \cap \nu = \alpha \} .$$
Then $\alpha ^\ast  _\nu < \nu$ because, by definition of $\beta (\nu )$, there exists a $\nu ^\prime < \nu$ and a $p \in J^X_{\rho (\nu )}$ such that $h_{\rho (\nu ),A_\nu}[\omega \times (J^X_{\nu ^\prime} \times \{ p \} )] \cap \nu$ is cofinal in $\nu$. But $p$ is in $h_{\rho (\nu ),A_\nu}[\omega \times (J^X_\nu \times \{ p_\nu \} )]$. So there is an $\alpha < \nu$ such that $h_{\rho (\nu ),A_\nu}[\omega \times (J^X_\alpha \times \{ p_\nu \} )] \cap \nu$ is cofinal in $\nu$. Thus $\alpha ^\ast _\nu  < \alpha < \nu$.
\smallskip \\
If $\nu \in S^+$, then we set $\alpha ^\ast _\nu := \alpha _\nu$.
\smallskip\\
For $\nu \in E$, let $f:\bar \nu \Rightarrow \nu$ iff, for some $f^\ast$,
\smallskip 

(1) $f=\langle \bar \nu ,f^\ast \upharpoonright J^D_{\mu _{\bar\nu}}, \nu \rangle$,
\smallskip 

(2) $f^\ast : I_{\mu _{\bar \nu}} \rightarrow I_{\mu _\nu}$ is $\Sigma _{n(\nu )}$-elementary,
\smallskip 

(3) $\alpha ^\ast_\nu$, $p_\nu$, $\alpha^\ast _{\mu _\nu}$, $P_\nu$ $\in rng(f^\ast )$,
\smallskip 

(4) $\nu \in rng(f^\ast )$ if $\nu < \mu _\nu$,
\smallskip 

(5) $f(\bar \nu )=\nu$ and $\bar \nu \in S^+ \Leftrightarrow \nu \in S^+$.
\medskip \\
By this, $\frak{F}$ is defined.
\medskip\\
Set $D=X$.
\medskip \\
Let $P^\ast _\nu$ be minimal such that $h^{n(\nu )-1}_{\mu _\nu}(i,P^\ast _\nu )=P_\nu$ for an $i \in \omega$.
\smallskip \\ 
Let $\alpha ^{\ast\ast} _{\mu _\nu}$ be minimal such that $h^{n(\nu )-1}_{\mu _\nu}(i,\alpha ^{\ast\ast} _{\mu _\nu} )=\alpha ^\ast _{\mu _\nu}$ for some $i \in \omega$.
\smallskip\\
Set
\smallskip

$\nu ^\ast = \emptyset$ if $\nu =\rho (\nu )$ 
\smallskip 

$\nu ^\ast = \nu$ if $\nu < \rho (\nu )$.
\medskip \\
For $\tau \in On$, let $S_\tau$ be defined as in lemma 10. For $\tau \in On$, $E_i \subseteq S_\tau$ and a $\Sigma _0$ formula $\varphi$, let
\smallskip \\
$h^\varphi _{\tau , E_i} (x_1, \dots ,x_m)$ the least $x_0 \in S_\tau$ w.r.t. the canonical well-ordering such that $\langle S_\tau , E_i \rangle \models \varphi (x_i)$ if such an element exists, 
\smallskip \\
and
\smallskip \\
$h^\varphi _{\tau , E_i} (x_1, \dots , x_m) = \emptyset$ else.
\smallskip \\
For $\tau \in On$ such that $\nu ^\ast ,\alpha ^\ast _\nu ,p_\nu ,\alpha ^{\ast\ast} _{\mu _\nu},P^\ast _\nu \in S_\tau$, let $H_\nu (\alpha ,\tau )$ be the closure of $S _\alpha \cup \{ \nu ^\ast ,\alpha ^\ast _\nu  ,p_\nu ,\alpha ^{\ast\ast} _{\mu _\nu},P^\ast _\nu \}$ under all $h^\varphi _{\tau , X \cap S_\tau ,  A_\nu \cap S_\tau }$. Then $H_\nu (\alpha ,\tau )\prec _1 \langle S_\tau , X \cap S_\tau , A_\nu \cap S_\tau , \{ \nu ^\ast ,\alpha ^\ast _\nu ,p_\nu ,\alpha ^{\ast\ast} _{\mu _\nu},P^\ast _\nu \} \rangle$ by the definition of $h^\varphi _{\tau , X \cap S_\tau ,  A_\nu \cap S_\tau }$. Let $M_\nu (\alpha ,\tau )$ be the collapse of $H_\nu (\alpha ,\tau )$. Let $\tau _0$ be the minimal $\tau$ such that  $\nu ^\ast ,\alpha ^\ast _\nu ,p_\nu , \alpha ^{\ast\ast} _{\mu _\nu},P^\ast _\nu \in S_\tau$. Define by induction for $\tau _0 \leq \tau < \rho (\nu )$:
\smallskip

$\alpha (\tau _0)=\alpha _\nu$
\smallskip

$\alpha (\tau +1)=sup(M_\nu (\alpha (\tau ),\tau +1) \cap \nu )$
\smallskip

$\alpha (\lambda )=sup\{ \alpha (\tau )\mid \tau < \lambda \}$ if $\lambda \in Lim$.
\smallskip \\
Set
\smallskip

$B_\nu = \{\langle \alpha (\tau ), M_\nu (\alpha (\tau ) ,\tau )\rangle \mid \tau _0 < \tau \in \rho (\nu )\}$ if $\nu < \rho (\nu )$,
\smallskip 

$B_\nu = \{ 0 \} \times A_\nu \cup \{ \langle 1,\nu ^\ast ,\alpha ^\ast _\nu ,p_\nu , \alpha ^{\ast\ast} _{\mu _\nu},P^\ast _\nu \rangle \}$ else.
\medskip\\
{\bf Lemma 22}
\smallskip \\
$B_\nu \subseteq J^X_\nu$ and $\langle I^0_\nu , B_\nu \rangle$ is rudimentary closed.
\smallskip \\
{\bf Proof:} If $\nu =\rho (\nu )$, then both claims are clear. Otherwise, we first prove $M^\nu (\alpha ,\tau ) \in J^X_\nu$ for all $\alpha < \nu$ and all $\tau \in \rho (\nu )$ such that $\tau _0 \leq \tau < \rho (\nu )$. Let such a $\tau$ be given and $\tau ^\prime \in \rho (\nu )-Lim$ be such that $X \cap S_\tau , A_\nu \cap S_\tau \in S_{\tau ^\prime}$ (rudimentary closedness of $\langle I^0_{\rho (\nu )}, A_\nu \rangle$). Let $\eta := sup(\tau ^\prime \cap Lim)$. Let $H$ be the closure of $\alpha \cup \{ \nu ^\ast ,\alpha ^\ast _\nu , p_\nu,\alpha ^{\ast \ast}_{\mu _\nu},P^\ast _\nu ,X \cap S_\tau,S_\tau ,A_\nu \cap S_\tau ,\eta \}$ under all $h^\varphi _{\tau ^\prime}$. Let $\sigma :H \cong S$ be the collapse of $H$ and $\sigma (\eta )=\bar \eta$. If $\eta \in S^X$, then $S = S_{\bar \tau ^\prime}$ for some $\bar \tau ^\prime$ by the condensation property of $L[X]$. If $\eta \not\in S^X$, then $S =S^{X \upharpoonright \bar \eta}_{\bar \tau ^\prime}$ for some $\bar \tau ^\prime$ where $S^{X \upharpoonright \bar \eta}_{\bar \tau ^\prime}$ is defined like $S_{\bar \tau ^\prime}$ with $X \upharpoonright \bar \eta$ instead of $X$. The reason is that, even if $\eta \notin S^X$, it is the supremum of points in $S^X$, because $S^X=\{ \beta (\nu ) \mid \nu$ singular in $L_\kappa [X]\}$.  In both cases, $S \in J^X_\nu$ and there is a function in $I_{\bar \eta + \omega}$ that maps $\alpha \cup \{\sigma ( \nu ^\ast ),\sigma (\alpha ^\ast _\nu ),\sigma ( p_\nu ),\sigma (\alpha ^{\ast \ast}_{\mu _\nu}),\sigma (P^\ast _\nu ),\sigma (X \cap S_\tau ),\sigma (S_\tau ),\sigma (A_\nu \cap S_\tau ),\sigma (\eta ) \}$ onto $S$. So $\nu$ would be  singular in $J^X_{\rho _\nu}$ if $\nu \leq \bar \tau ^\prime$. But this contradicts the definition of $\beta (\nu )$. Therefore, $\sigma ( \nu ^\ast ),\sigma (\alpha ^\ast _\nu ),\sigma ( p_\nu ),\sigma (\alpha ^{\ast \ast}_{\mu _\nu}),\sigma (P^\ast _\nu ),\sigma (X \cap S_\tau ),\sigma (S_\tau ),\sigma (A_\nu \cap S_\tau ),\sigma (\eta )\in J^X_\nu$. Let $\bar H _\nu (\alpha ,\tau )$ be the closure of $S_{\alpha} \cup \{ \sigma ( \nu ^\ast ),\sigma (\alpha ^\ast _\nu ),\sigma ( p_\nu ),$ $\sigma (\alpha ^{\ast \ast}_{\mu _\nu}),\sigma (P^\ast _\nu ),\sigma (X \cap S_\tau ),\sigma (S_\tau ),\sigma (A_\nu \cap S_\tau ),\sigma (\eta ) \}$ under all $h^\varphi _{\sigma (S_\tau ), \sigma (X \cap S_\tau ),\sigma (A_\nu \cap S_\tau )}$ where these are defined like $h^\varphi _{\tau ,E_i}$ but with $\sigma (S_\tau )$ instead of $S_\tau$. Then $\bar H_\nu (\alpha ,\tau ) \prec _1 \langle \sigma (S_\tau ), \sigma (X \cap S_\tau ),\sigma (A_\nu \cap S_\tau ), \{ \sigma ( \nu ^\ast ),\sigma (\alpha ^\ast _\nu ),\sigma ( p_\nu ),$ $\sigma (\alpha ^{\ast \ast}_{\mu _\nu}),\sigma (P^\ast _\nu ),\sigma (X \cap S_\tau ),$ $\sigma (S_\tau ),\sigma (A_\nu \cap S_\tau ),\sigma (\eta ) \} \rangle$ and $M_\nu (\alpha ,\tau )$ is the collapse of $\bar H_\nu (\alpha ,\tau )$. Since $\nu < \rho (\nu )$ and $\nu$ is a cardinal in $I_{\beta (\nu )}$, $J^X_\nu \models ZF^-$. So we can form the collapse inside $J^X_\nu$. Thus $M_\nu (\alpha ,\tau ) \in J^X_\nu$.

Now, we turn to rudimentary closedness. Since $B_\nu$ is unbounded in $\nu$, it suffices to prove that the initial segments of $B_\nu$ are elements of $J^X _\nu$. Such an initial segment is of the form $\langle M_\nu (\alpha (\tau ),\tau ) \mid \tau < \gamma \rangle$ where $\gamma < \rho (\nu )$, and we have $H_\nu (\alpha (\tau ),\delta _\tau ) = H_\nu (\alpha (\tau ) ,\tau )$ where $\delta _\tau$ is for $\tau < \gamma$ the least $\eta \geq \tau$ such that $\eta \in H_\nu (\alpha (\tau ),\gamma ) \cup \{ \gamma \}$. Since $\delta _\tau \in H_\nu (\alpha (\tau ), \gamma ) \prec _1 \langle S_\gamma , X \cap S_\gamma , A_\nu \cap S_\gamma , \{ \dots \} \rangle$, $(H_\nu (\alpha (\tau ),{\delta _\tau })) ^{H_\nu (\alpha (\gamma ) ,\gamma)} = H_\nu (\alpha (\tau ),\tau )$. Let $\pi : M_\nu (\alpha (\gamma ),\gamma ) \rightarrow S_\gamma$ be the uncollapse of $H_\nu (\alpha (\gamma ),\gamma )$. Then, by the $\Sigma _1$-elementarity of $\pi$, $M_\nu (\alpha (\tau ),\tau ) =M_\nu (\alpha (\tau ) ,{\delta _\tau})$ is the collapse of  $(H(\alpha (\tau ),{\pi ^{-1}(\delta _\tau )}) )^{M_\nu (\alpha (\gamma ),\gamma)}$. So $\langle M_\nu (\alpha (\tau ),\tau ) \mid \tau < \gamma \rangle$ is definable from $M_\nu (\alpha (\gamma ),\gamma ) \in J^X_\nu$. $\Box$
\bigskip\\
{\bf Lemma 23}
\smallskip \\
For $x,y_i \in J^X_{\nu}$, the following are equivalent:
\smallskip \\
(i) $x$ is $\Sigma _1$-definable in $\langle I^0_{\rho (\nu )},A_\nu \rangle$ with the parameters $y_i,\nu ^\ast ,\alpha ^\ast _\nu ,p_\nu , \alpha ^{\ast\ast} _{\mu _\nu},P^\ast _\nu$.
\smallskip \\
(ii) $x$ is $\Sigma _1$-definable in $\langle I^0_\nu , B_\nu \rangle$ with the parameters $y_i$.
\smallskip \\
{\bf Proof:} For $\nu =\rho (\nu )$, this is clear. Otherwise, let first $x$ be uniquely determined in $\langle I^0_{\rho (\nu )} , A_\nu \rangle$ by $(\exists z) \psi (z,x, \langle y_i , \nu ^\ast ,\alpha ^\ast _\nu ,p_\nu , \alpha ^{\ast\ast} _{\mu _\nu},P^\ast _\nu \rangle )$ where $\psi$ is a $\Sigma _0$ formula. That is equivalent to $(\exists \tau)(\exists z \in S_\tau) \psi (z,x,\langle y_i , \nu ^\ast ,\alpha ^\ast _\nu ,p_\nu , \alpha ^{\ast\ast} _{\mu _\nu},P^\ast _\nu \rangle )$ and that again to $(\exists \tau )H_\nu (\alpha (\tau ), \tau ) \models (\exists z) \psi (z,x, \langle y_i , \nu ^\ast ,\alpha ^\ast _\nu ,p_\nu , \alpha ^{\ast\ast} _{\mu _\nu},P^\ast _\nu \rangle )$. If $\tau$ is large enough, the $y_i$ are not moved by the collapsing map, since then $y_i \in J^X_{\alpha (\tau )} \subseteq H_\nu (\alpha (\tau ),\tau )$. Let $\bar \nu ,\alpha ,p,\alpha ^\prime ,P$ be the images of $\nu ^\ast ,\alpha ^\ast _\nu ,p_\nu , \alpha ^{\ast\ast} _{\mu _\nu},P^\ast _\nu$ under the collapse. Then $(\exists \tau )(y_i \in J^X_{\alpha (\tau )}$ $and$ $M_\nu (\alpha (\tau ),\tau ) \models (\exists z) \psi (z,x, \langle y_i ,\bar \nu ,\alpha ,p,\alpha ^\prime ,P\rangle ))$ defines $x$. So it is definable in $\langle I^0,B_\nu \rangle$.\\
Since $B_\nu$ and the satisfaction relation of $\langle I^0_\gamma , B \rangle$ are $\Sigma _1$-definable over $\langle I^0_{\rho (\nu )}, A_\nu \rangle$, the converse is clear. $\Box$
\medskip \\
{\bf Lemma 24}
\smallskip \\
Let $H \prec _1 \langle I^0_\nu , B_\nu \rangle$ for a $\nu \in E$ and $\pi : \langle I^0_\mu , B \rangle \rightarrow \langle I^0_\nu , B_\nu \rangle$ be the uncollapse of $H$. Then $\mu \in E$ and $B=B_\mu$.
\smallskip \\
{\bf Proof:} First, we extend $\pi$ like in lemma 19. Let 
\smallskip \\
$M$ $=\{ x \in J^X _{\rho (\nu )} \mid x$ is $\Sigma _1$-definable in $\langle I^0_{\rho (\nu )},A_\nu \rangle$ with parameters from $rng(\pi ) \cup \{ p_\nu  ,\nu ^\ast ,\alpha ^\ast _\nu ,\alpha ^{\ast \ast}_{\mu _\nu},P^\ast _\nu\}$ $\}$.
\smallskip \\
Then $rng(\pi ) = M \cap J^X_\nu$. For, if $x \in M \cap J^X_\nu$, then there are by definition of $M$ $y_i \in rng(\pi )$ such that $x$ is $\Sigma _1$-definable in $\langle I^0_{\rho (\nu )},A_\nu \rangle$ with the parameters $y_i$ and $p_\nu,\nu ^\ast ,\alpha ^\ast _\nu ,\alpha ^{\ast \ast}_{\mu _\nu},P^\ast _\nu$. Thus it is $\Sigma _1$-definable in $\langle I^0_\nu , B_\nu \rangle$ with the $y_i$ by lemma 23. Therefore, $x \in rng(\pi )$ because $y_i \in rng(\pi ) \prec _1 \langle I^0_\nu , B_\nu \rangle$. Let $\hat \pi : \langle I^0_\rho , A \rangle \rightarrow \langle I^0_{\rho (\nu )},A_\nu \rangle$ be the uncollapse of $M$. Then $\hat \pi$ is an extension of $\pi$, since $M \cap J^X_\nu$ is an $\in$-initial segment of $M$ and $rng(\pi ) = M \cap J^X _\nu$. In addition, there is by lemma 19 a $\Sigma _{n(\nu )}$- elementary extension $\tilde \pi : I_\beta \rightarrow I_{\beta (\nu )}$ such that $\rho$ is the $(n(\nu )-1)$-th projectum of $I_\beta$ and $A$ is the $(n(\nu )-1)$-th standard code of it. Let $\tilde \pi (p) = p_\nu$ and $\tilde \pi (\alpha )=\alpha ^\ast _\nu$. And we have $\tilde \pi (\mu )=\nu$ if $\nu < \beta (\nu )$. In this case, $\nu \in rng(\pi )$ by the definition of $\nu ^\ast$. Since $\tilde\pi$ is $\Sigma _1$-elementary, cardinals of $J^X_\mu$ are mapped on cardinals of $J^X_\nu$.

 Assume $\nu \in S^+$. Suppose there was a cardinal $\tau > \alpha$ of $J^X_\mu$. Then $\pi (\tau )>\alpha _\tau$ was a cardinal in $J^X_\nu$. But this is a contradiction.

Next, we note that $\mu$ is $\Sigma _{n(\nu )}$-singular over $I_\beta$. If $\nu \in S^+$, then, by the definition of $p_\nu$, $J^X_\rho = h_{\rho,A}[\omega \times (\alpha \times \{ p \} )]$ is clear. So there is an over $\langle I^0_\rho , A \rangle$ $\Sigma _1$-definable function from $\alpha$ cofinal into $\mu$. But since $\rho$ is the $(n(\nu )-1)$-th projectum and $A$ is the $(n(\nu )-1)$-th code of it, this function is $\Sigma _n$-definable over $I_\beta$. Now, suppose $\nu \notin S^+$. Let $\lambda := sup (\pi [\mu ])$. Since $\lambda > \alpha ^\ast _\nu$, there is a $\gamma < \lambda$ such that
$$sup(h_{\rho (\nu ),A_\nu}[\omega \times (J^X_\gamma \times \{ q_\nu \} )] \cap \nu ) \geq \lambda .$$
And since $rng(\pi )$ is cofinal in $\lambda$, there is such a $\gamma \in rng(\pi )$. Let $\gamma =\pi (\bar \gamma )$. By the $\Sigma _1$-elementarity of $\tilde \pi$, $\bar \gamma < \mu$ and setting $\tilde \pi (q)=q_\nu$ we have for every $\eta < \mu$
$$\langle I_\rho ,A \rangle \models (\exists x \in J^X_{\bar \gamma})(\exists i)h_{\rho ,A}(i,\langle x,p \rangle )>\eta . $$ 
Hence $h_{\rho ,A}[\omega \times (J^X_{\bar \gamma} \times \{ q \} )]$ is cofinal in $\mu$. This shows $\mu \in E$.

On the other hand, $\mu$ is $\Sigma _{n(\nu )-1}$-regular over $I_\beta$ if $n(\nu )>1$. Assume there was an over $I_\beta$ $\Sigma _{n(\nu )-1}$-definable function $f$ and some $x \in \mu$ such that $f[x]$ was cofinal in $\mu$. I.e. $(\forall y \in \mu )(\exists z \in x )(f(x)>y)$ would hold in $I_\beta$. Over $I_\beta$, $(\exists z \in x)(f(z)>y)$ is $\Sigma _{n(\nu )-1}$. So it is $\Sigma _0$ over $\langle I^0_\rho , A \rangle$. But then also $(\forall y \in \mu )(\exists z \in x)(f(z)>y)$ is $\Sigma _0$ over $\langle I^0_\rho , A \rangle$ if $\mu < \rho$. Hence it is $\Sigma _{n(\nu )}$ over $I_\beta$. But then the same would hold for $\tilde \pi (x)$ in $I_{\beta (\nu )}$. This contradicts the definition of $n(\nu )$! Now, let $\mu = \rho$. Since $\alpha$ is the largest cardinal in $I_\mu$, we had in $f$ also an over $I_\beta$ $\Sigma _{n(\nu )-1}$-definable function from $\alpha$ onto $\rho$ and therefore one from $\alpha$ onto $\beta$. But this contradicts lemma 21 and the fact that $\rho$ is the $(n(\nu )-1)$-th projectum of $\beta$. If $n(\nu )=1$, then we get with the same argument that $\mu$ is regular in $I_\beta$.

The previous two paragraphs show $\beta = \beta (\mu )$ and $n(\mu )=n(\nu )$. We are done if we can also show that $\alpha =\alpha ^\ast _\mu, \pi (\alpha ^{\ast\ast}_{\mu _\mu})=\alpha ^{\ast\ast}_{\mu _\nu}, p=p_\mu , \pi (P^\ast _\mu )=P^\ast _\nu$, because $\tilde \pi$ is $\Sigma _1$-elementary, $\tilde \pi (h^\varphi _{\tau , X \cap S_\tau , A_\mu \cap S_\tau} (x_i))=h^\varphi _{\tilde \pi (\tau ) , X \cap S_{\tilde \pi (\tau )} , A_\nu \cap S_{\tilde \pi (\tau )}} (x_i)$ for all $\Sigma _1$ formulas $\varphi$ and $x_i \in S_\tau$.

For $\nu \in S^+$, $\alpha =\alpha _\mu$ was shown above. So let $\nu \notin S^+$. By the $\Sigma _1$-elementarity of $\tilde \pi$, we have for all $\alpha \in \mu$
$$h_{\rho ,A}[\omega \times (J^X_\alpha \times \{ p \} )] \cap \mu = \alpha \Leftrightarrow
h_{\rho (\nu ),A_\nu}[\omega \times (J^X_{\pi (\alpha )} \times \{ p_\nu \} )] \cap \nu = \pi (\alpha ).$$ 
The same argument proves  $\pi (\alpha ^{\ast\ast}_{\mu _\mu})=\alpha ^{\ast\ast}_{\mu _\nu}$. Finally, $p=p_\mu$ and $\pi (P^\ast _\mu )=P^\ast _\nu$ can be shown as in (5) in the proof of lemma 19. $\Box$  
\medskip \\
{\bf Lemma 25}
\smallskip \\
Let $H \prec _1 \langle I^0_\nu , B_\nu \rangle$ and $\lambda = sup(H \cap \nu )$ for a $\nu \in E$. Then $\lambda \in E$ and $B_\nu \cap J^X_\lambda = B_\lambda$.
\smallskip \\
{\bf Proof:} Let $\pi _0 : \langle I^0_\mu , B_\mu \rangle \rightarrow \langle I^0_\lambda , B_\nu \cap J^X_\lambda \rangle$ be the uncollapse of $H$ and let $\pi _1 : \langle I^0_\lambda , B_\nu \cap J^X_\lambda \rangle \rightarrow \langle I^0_\nu , B_\nu \rangle$ be the identity. Since $L[X]$ has coherence, $\pi _0$ and $\pi _1$ are $\Sigma _0$-elementary. By lemma 18, $\pi _0$ is even $\Sigma _1$-elementary, because it is cofinal. To show $B_\lambda =B_\nu \cap J^X_\lambda$, we extend $\pi _0$ and $\pi _1$ to $\hat \pi _0 : \langle I^0_{\rho (\mu )}, A_\mu \rangle \rightarrow \langle I^0_\rho , A \rangle$ and $\hat \pi _1 : \langle I^0_\rho , A \rangle \rightarrow \langle I^0_{\rho (\nu )},A_\nu \rangle$ in such a way that $\hat \pi _0$ is $\Sigma _1$-elementary and $\hat \pi _1$ is $\Sigma _0$-elementary. Then we know from lemma 19 that $\rho$ is the $(n(\nu )-1)$-th projectum  of some $\beta$ and $A$ is the $(n(\nu )-1)$-th code of it. So there is a $\Sigma _{n(\nu )}$-elementary extension of $\tilde \pi _0 : I_{\bar\beta} \rightarrow I_\beta$. We can again use the argument from lemma 24 to show that $\lambda$ is $\Sigma _{n(\nu )-1}$-regular over $I_\beta$. But on the other hand, $\lambda$ is as supremum of $H \cap On$ $\Sigma _{n(\nu )}$-singular over $I_\beta$. From this, we conclude as in the proof of lemma 24 that $B_\lambda =B_\nu \cap J^X_\lambda$.

First, suppose $\nu \in S^+$. Since $\alpha _\nu \in H \prec _1 \langle I^0_\nu , B_\nu \rangle$, $\alpha _\nu < \lambda \leq \nu$. Since $I_\nu \models (\alpha _\nu$ is the largest cardinal), we therefore have $\lambda \notin Card$. In addition, $\alpha _\nu$ is the largest cardinal in $I_\lambda$. Assume $\tau$ was the next larger cardinal. Then $\tau$ was $\Sigma _1$-definable in $I_\lambda$ with parameter $\alpha _\nu$ and some $\tau ^\prime \in H$ and hence it was in $H$. By the $\Sigma _1$-elementarity of $\pi _0$, $\pi _0^{-1}(\tau ) > \pi _0^{-1}(\alpha _\nu )
= \alpha _\mu$ was also a cardinal in $I_\mu$. But this contradicts the definition of $\alpha _\mu$.

But now to $B_\lambda = B_\nu \cap J^X_\lambda$. First, assume $\nu \notin S^+$. Let $\pi = \pi _1 \circ \pi _0 : \langle I^0_\mu , B _\mu \rangle \rightarrow \langle I^0_\nu , B_\nu \rangle$ and $\hat \pi : \langle I^0_{\rho (\mu )} , A_\mu \rangle \rightarrow \langle I^0_{\rho (\nu )},A_\nu \rangle$ be the extension constructed in the proof of lemma 24. Let $\gamma = sup(rng( \hat \pi ))$. Then $\hat \pi ^\prime = \hat \pi \cap (J^X_{\rho (\mu )} \times J^X_\gamma ) : \langle I^0_{\rho (\mu )},A_\mu \rangle \rightarrow \langle I^0_\gamma , A_\nu \cap J^X_\gamma \rangle$ is $\Sigma _0$-elementary, by coherence of $L_\kappa [X]$, and cofinal. Thus $\hat \pi ^\prime$ is $\Sigma _1$-elementary. Let $H ^\prime = h_{\gamma , A_\nu \cap J^X_\gamma}[ \omega \times ( J^X_\lambda \times \{ p_\nu \} ) ]$ and $\hat \pi _1 : \langle I^0_\rho , A \rangle \rightarrow \langle I^0_{\rho (\nu ) },A_\nu \rangle$ be the uncollapse of $H^\prime$. Then $H= rng( \hat \pi ^\prime ) \subseteq H^\prime$. To see this, let $z \in rng(\hat \pi ^\prime )$ and $z=\hat \pi ^\prime (y)$. Then by definition of $p_\mu$, there is an $x \in J^X_\mu$ and an $i \in \omega$ such that $y=h_{\rho (\mu ),A_\mu}(i,\langle x,p_\mu \rangle )$. By the $\Sigma _1$-elementarity of $\hat \pi ^\prime$, we therefore have $z=h_{\gamma , A_\nu \cap J^X_\gamma}(i,\langle  \hat \pi ^\prime (x), \hat \pi ^\prime (p_\mu ) \rangle )$. But $\hat \pi ^\prime (p_\mu )=\hat \pi (p_\mu ) =p_\nu$ and $\hat \pi ^\prime (x) \in J^X_\lambda$.

In addition, $sup(H^\prime \cap \nu ) = \lambda$. That $\sup (H^\prime \cap \nu ) \geq \lambda$ is clear. Conversely, let $x \in H^\prime \cap \nu$, i.e. $x=h_{\gamma , A_\nu \cap J^X_\gamma}(i, \langle y,p_\nu \rangle )$ for some $i \in \omega$ and a $y \in J^X_\lambda$. Then $x$ is uniquely determined by $\langle I^0_\gamma , A_\nu \cap J^X_\gamma \rangle \models (\exists z) \psi _i (z,x, \langle y,p_\nu \rangle )$. But such a $z$ exists already in a $H^0_\nu (\alpha , \tau )$ where $H^0_\nu (\alpha ,\tau )$ is the closure of $S _\alpha$ under all $h^\varphi _{\tau , X \cap S_\tau ,  A_\nu \cap S_\tau }$. Since $\gamma = sup(rng(\hat \pi ))$ and $\lambda = sup(rng(\pi ))$ we can pick such $\tau \in rng(\hat \pi )$ and $\alpha \in rng(\pi )$. Let $\bar \tau = \hat \pi ^{-1}(\tau )$ and  $\bar \alpha = \hat \pi ^{-1}(\alpha )$. Let $\vartheta = sup(\nu \cap H^0_\nu (\alpha ,\tau ))$ and $\bar \vartheta = sup (\mu \cap H^0 _\mu (\bar \alpha ,\bar \tau ))$. Since $\nu$ is regular in $I_{\rho (\nu )}$, $\vartheta  < \nu$. Analogously, $\bar \vartheta < \mu$. But of course $\hat \pi (\bar \vartheta )= \vartheta$. So $x < \vartheta =\hat \pi (\bar \vartheta ) < sup(\hat \pi [\mu ])=\lambda$. 

If $\nu \in S^+$, we may define $H^\prime$ as $ h_{\gamma , A_\nu \cap J^X_\gamma}[ \omega \times ( J^X_{\alpha _\nu} \times \{ p_\nu \} ) ]$ and still conclude that $H= rng( \hat \pi ^\prime ) \subseteq H^\prime$ and $sup(H^\prime \cap \nu ) = \lambda$ by the definition of $p_\nu$.

By lemma 19, $\hat \pi : \langle I^0_\rho , A \rangle \rightarrow \langle I^0_{\rho (\nu )},A_\nu \rangle$ may be extended to a $\Sigma _{n(\nu )-1}$-elementary embedding $\tilde \pi _1 : I_\beta \rightarrow I_{\beta (\nu )}$ such that $\rho$ is the $(n(\nu ) -1)$-th projectum of $I_\beta$ and $A$ is the $(n(\nu )-1)$-th standard code of it. Let $\hat \pi _0 = \hat \pi ^{-1}_1 \circ \hat \pi$. Then $\hat \pi _0 : \langle I^0_{\rho (\mu )},A_\mu \rangle \rightarrow \langle I^0_\rho , A \rangle$ is $\Sigma _0$-elementary, by the coherence of $L_\kappa [X]$, and cofinal. Thus it is $\Sigma _1$-elementary by lemma 18. Applying again lemma 19, we get a $\Sigma _{n(\nu )}$-elementary $\tilde \pi _0 : I_{\beta (\mu )} \rightarrow I_\beta$.

As in lemma 24, it suffices to prove $\beta = \beta (\lambda )$, $n(\nu )=n(\lambda )$, $\rho = \rho (\lambda )$, $A=A_\lambda$, $\hat \pi _1 ^{-1}(p_\nu )=p_\lambda$, $\hat \pi _1 ^{-1}(P^\ast _\nu )=P^\ast _\lambda$, $\alpha ^\ast _\nu = \alpha ^\ast _\lambda$ and $\hat \pi _1 ^{-1}(\alpha ^{\ast\ast}_{\mu _\nu} )=\alpha ^{\ast \ast}_{\mu _\lambda}$. So, if $n(\nu )>1$, we have to show that $\lambda$ is $\Sigma _{n(\nu )-1}$-regular over $I_\beta$. If $n(\nu )=1$, then $I_\beta \models (\lambda$ regular) suffices. In addition, $\lambda$ must be $\Sigma _{n(\nu )}$-singular over $I_\beta$. For regularity, consider $\tilde \pi _0$ and, as in lemma 24, the least $x \in \lambda$ proving the opposite if such an $x$ exists. This is again $\Sigma _n$-definable and therefore in $rng(\tilde \pi _0)$. But then $\tilde \pi _0^{-1} (x)$ had the same property in $I_{\beta (\mu )}$. Contradiction!  

Now, assume $\nu \in S^+$. Since $I_\nu \models (\alpha _\nu$ is the largest cardinal), $H^\prime \cap \nu$ is transitive. Thus $H^\prime \cap \nu = \lambda$. Since $\hat \pi _1 : \langle I^0_\rho , A \rangle \rightarrow \langle I^0_\gamma , A \cap J^X_\gamma \rangle$ is $\Sigma _1$-elementary and $\lambda \subseteq H^\prime = rng(\hat \pi _1)$, we have $\lambda =\lambda \cap h_{\rho , A}[\omega \times (J^X_{\alpha _\nu} \times \{ \hat \pi _1^{-1}(p_\nu ) \} ) ]$. I.e. there is a  $\Sigma _1$-map over $\langle I_\rho , A \rangle$ from $\alpha _\nu$ onto $\lambda$. But this is then $\Sigma _{n(\nu )}$-definable over $I_\beta$ and $\lambda$ is $\Sigma _{n(\nu )}$-singular over $I_\beta$.

If $\nu \notin S^+$, then the fact that $\lambda$ is $\Sigma _{n (\nu )}$-singular over $I_\beta$, $\alpha ^\ast _\nu = \alpha ^\ast _\lambda$ and  $\hat \pi _1 ^{-1}(\alpha ^{\ast\ast}_{\mu _\nu} )=\alpha ^{\ast \ast}_{\mu _\lambda}$ may be seen as in lemma 24 because $\pi _0 (\alpha ^\ast _\mu )=\alpha ^\ast _\nu \in rng(\pi _0)$.

That $\hat \pi _1^{-1} (p_\nu )=p_\lambda$ and  $\hat \pi _1 ^{-1}(P^\ast _\nu )=P^\ast _\lambda$ can again be proved as in (5) in the proof of lemma 19. $\Box$
\medskip \\
{\bf Lemma 26}
\smallskip \\
Let $\nu \in E$ and $\Lambda (\xi ,\nu )=\{ sup(h_{\nu ,B_\nu}[\omega \times (J^X_\beta \times \{ \xi \} )] \cap \nu) < \nu \mid \beta \in Lim \cap \nu \}$. Let $\bar \eta < \bar \nu$ and $\pi : \langle I^0_{\bar\nu} ,B \rangle \rightarrow \langle I^0_\nu ,B_\nu \rangle$ be $\Sigma _1$-elementary. Then $\Lambda (\bar \xi ,\bar \nu )\cap \bar \eta \in J^X_{\bar \nu}$ and $\pi (\Lambda (\bar \xi ,\bar \nu )\cap \bar \eta )=\Lambda (\xi ,\nu ) \cap \pi (\bar \eta )$ where $\pi (\bar \xi )=\xi$ and $\pi (\bar \eta )=\eta$.
\smallskip \\
{\bf Proof:} 
\smallskip \\
(1) Let $\lambda \in \Lambda (\xi ,\nu )$. Then $\Lambda (\xi ,\lambda )=\Lambda (\xi ,\nu )\cap \lambda$.
\smallskip \\
Let $\beta _0$ be minimal such that
\smallskip

$sup(h_{\nu ,B_\nu}[\omega \times (J^X_{\beta _0} \times \{ \xi \} )] \cap \nu )=\lambda$.
\smallskip \\
Then, by lemma 25, for all $\beta \leq \beta _0$
\smallskip 

$h_{\lambda ,B_\lambda}[\omega \times (J^X_\beta \times \{ \xi \} )] =
 h_{\nu,B_\nu}[\omega \times (J^X_\beta \times \{ \xi \} )]$
\smallskip \\
and for all $\beta _0 \leq \beta$
\smallskip 

$h_{\lambda ,B_\lambda}[\omega \times (J^X_{\beta _0} \times \{ \xi \} )] \subseteq
h_{\lambda ,B_\lambda}[\omega \times (J^X_\beta \times \{ \xi \} )]$
\smallskip

$\subseteq h_{\nu,B_\nu}[\omega \times (J^X_\beta \times \{ \xi \} )]$.
\smallskip \\
So $\Lambda (\xi ,\lambda )=\Lambda (\xi ,\nu )\cap \lambda$.
\smallskip \\
(2) $\Lambda (\bar \xi ,\bar \nu )\cap \bar \eta \in J^X_{\bar \nu}$
\smallskip \\
Let $\bar \lambda := sup(\Lambda (\bar \xi ,\bar \nu )\cap \bar \eta +1)$. Then, by (1), 
$\Lambda (\bar\xi ,\bar\nu )\cap \bar \eta +1 =\Lambda (\bar\xi ,\bar\lambda ) \cup \{ \bar \lambda \}$. But $\Lambda (\bar\xi ,\bar\lambda )$ is definable  over $I_{\beta (\bar \lambda )}$. Since $\beta (\bar \lambda )< \bar \nu$, we get $\Lambda (\bar \xi ,\bar \nu )\cap \bar \eta +1\in J^X_{\bar \nu}$.
\smallskip\\
(3) Let $sup(h_{\bar \nu ,B_{\bar \nu}}[\omega \times (J^X_{\bar\beta} \times \{ \bar\xi \} )]) < \bar\nu$ and $\pi (\bar \beta )=\beta$. Then 
\smallskip

$\pi (sup(h_{\bar \nu ,B_{\bar \nu}}[\omega \times (J^X_{\bar\beta} \times \{ \bar\xi \} )] \cap \bar \nu ))=sup(h_{\nu,B_\nu}[\omega \times (J^X_\beta \times \{ \xi \} )] \cap \nu )$.
\smallskip \\
Let $\bar \lambda :=sup(h_{\bar \nu ,B_{\bar \nu}}[\omega \times (J^X_{\bar\beta} \times \{ \bar\xi \} )] \cap \bar \nu )$. Then $\langle I^0_{\bar \nu},B_{\bar \nu}\rangle \models \neg (\exists \bar \lambda < \theta )(\exists i \in \omega )(\exists \xi _i < \bar \beta )(\theta = h_{\bar \nu ,B_{\bar \nu}}(i,\langle \xi _i,\bar \xi \rangle ))$. So $\langle I^0_\nu ,B_\nu \rangle \models \neg (\exists \lambda < \theta )(\exists i \in \omega )(\exists \xi _i < \beta )(\theta = h_{\nu ,B_\nu}(i,\langle \xi _i,\xi \rangle ))$ where $\pi (\bar \lambda )=\lambda$. I.e. $sup(h_{\nu ,B_\nu}[\omega \times (J^X_\beta \times \{ \xi \} )] \cap \nu ) \leq \lambda$. But $(\pi \upharpoonright J^X_{\bar\lambda}): \langle I^0_{\bar \lambda},B_{\bar \lambda}\rangle \rightarrow \langle I^0_\lambda ,B_\lambda \rangle$ is elementary. So, if $\langle I^0_{\bar \lambda},B_{\bar \lambda}\rangle \models (\forall \eta )(\exists \xi _i \in \bar \beta )(\exists n \in \omega )(\eta \leq h_{\bar \lambda ,B_{\bar \lambda}}(n,\langle \xi _i,\bar \xi \rangle ))$, then $\langle I^0_{\lambda},B_{\lambda}\rangle \models (\forall \eta )(\exists \xi _i \in \beta )(\exists n \in \omega )(\eta \leq h_{\lambda ,B_{\lambda}}(n,\langle \xi _i,\xi \rangle ))$. But by lemma 25, $h_{\lambda ,B_\lambda}[\omega \times (J^X_\beta \times \{ \xi \} )] \subseteq h_{\nu,B_\nu}[\omega \times (J^X_\beta \times \{ \xi \} )]$. I.e. it is indeed $\lambda = sup(h_{\nu,B_\nu}[\omega \times (J^X_\beta \times \{ \xi \} )] \cap \nu )$.
\smallskip \\
(4) $\pi (\Lambda (\bar \xi ,\bar \nu )\cap \bar \eta )=\Lambda (\xi ,\nu ) \cap \pi (\bar \eta )$
\smallskip \\
For $\bar \lambda \in \Lambda (\bar \xi ,\bar \nu )$,
\smallskip \\
$\pi (\Lambda (\bar \xi ,\bar \nu ) \cap \bar \lambda )$
\smallskip

by (1)
\smallskip \\
$=\pi (\Lambda (\bar \xi , \bar \lambda ))$
\smallskip 

by $\Sigma _1$-elementarity of $\pi$
\smallskip \\
$=\Lambda (\xi ,\pi (\bar \lambda ) )$
\smallskip

by (1) and (3)
\smallskip \\
$=\Lambda (\xi ,\nu ) \cap \pi (\bar \lambda )$.
\smallskip \\
So, if $\Lambda (\bar \xi ,\bar \nu )$ is cofinal in $\bar \nu$, then we are finished. But if there exists $\bar \lambda :=max(\Lambda (\bar\xi ,\bar\nu ))$, then, by (1) and (2), $\Lambda (\bar\xi ,\bar\nu ) \in J^X_{\bar \nu}$, and it suffices to show $\pi (\Lambda (\bar \xi ,\bar \nu ))=\Lambda (\xi ,\nu )$. To this end, let $\bar \beta$ be maximal such that $\bar \lambda =sup(h_{\bar \nu ,B_{\bar \nu}}[\omega \times (J^X_{\bar\beta} \times \{ \bar\xi \} )] \cap \bar \nu )$. I.e. 
$h_{\bar \nu ,B_{\bar \nu}}[\omega \times (J^X_{\bar\beta +1} \times \{ \bar\xi \} )]$ is cofinal in $\bar \nu$. So, since $\pi [h_{\bar \nu ,B_{\bar \nu}}[\omega \times (J^X_{\bar\beta +1} \times \{ \bar\xi \} )]] \subseteq h_{\nu,B_\nu}[\omega \times (J^X_{\beta +1} \times \{ \xi \} )]$ where $\pi (\bar \beta )=\beta$,  $sup(rng(\pi )\cap \nu ) \leq sup(h_{\nu,B_\nu}[\omega \times (J^X_{\beta +1} \times \{ \xi \} )] \cap \nu)$. Hence indeed $\pi (\Lambda (\bar \xi ,\bar \nu ))=\Lambda (\xi ,\nu )$. $\Box$
\medskip \\
{\bf Lemma 27}
\smallskip \\
Let $\nu \in E$, $H \prec _1 \langle I^0_\nu ,B_\nu \rangle$ and $\lambda = sup(H \cap \nu )$. Let $h: I^0_{\bar \lambda} \rightarrow I^0_\lambda$ be $\Sigma _1$-elementary and $H \subseteq rng(h)$. Then $\lambda \in E$ and $h:\langle I^0_{\bar \lambda},B_{\bar \lambda}\rangle \rightarrow \langle I^0_\lambda ,B_\lambda \rangle$ is $\Sigma _1$-elementary.
\smallskip \\
{\bf Proof:} By lemma 25, $B_\lambda =B_\nu \cap J^X_\lambda$. So it suffices, by lemma 24, to show $rng(h) \prec _1 \langle I^0_\lambda ,B_\lambda \rangle$. Let $x_i \in rng(h)$ and $\langle I^0_\lambda ,B_\lambda \rangle \models (\exists z)\psi (z,x_i)$ for a $\Sigma _0$ formula $\psi$. Then we have to prove that there exists a $z \in rng(h)$ such that $\langle I^0_\lambda ,B_\lambda \rangle \models \psi (z,x_i)$. Since $\lambda =sup(H \cap \nu )$, there is a $\eta \in H \cap Lim$ such that $\langle I^0_\eta ,B_\lambda \cap J^X_\eta \rangle \models (\exists z) \psi (z,x_i)$. And since $H \prec _1 \langle I^0_\nu ,B_\nu \rangle$, we have $\langle I^0_\eta ,B_\lambda \cap J^X_\eta \rangle \in H \subseteq rng(h)$. So also
$$rng(h) \models (\langle I^0_\eta ,B_\lambda \cap J^X_\eta \rangle \models (\exists z) \psi (z,x_i))$$
because $rng(h) \prec _1 I^0_\lambda$.
Hence there is a $z \in rng(h)$ such that $\langle I^0_\eta ,B_\lambda \cap J^X_\eta \rangle \models \psi (z,x_i)$. I.e. $\langle I^0_\lambda ,B_\lambda \rangle \models \psi (z,x_i)$. $\Box$ 
\bigskip \\
{\bf Lemma 28}
\smallskip \\
Let $f:\bar \nu \Rightarrow \nu$, $\bar \nu \sqsubset \bar \tau \sqsubseteq \mu _{\bar \nu}$ and $f(\bar \tau )=\tau$. If $\bar \tau \in S^+ \cup \widehat{S}$ is independent, then $(f \upharpoonright J^D_{\alpha _{\bar \tau}}):\langle J^D_{\alpha _{\bar \tau}},D_{\alpha _{\bar \tau}},K_{\bar \tau} \rangle \rightarrow \langle J^D_{\alpha _{\tau}},D_{\alpha _{\tau}},K_\tau \rangle$ is $\Sigma _1$-elementary.
\smallskip \\
{\bf Proof:} If $\bar \tau =\mu _{\bar \tau} < \mu _{\bar \nu}$, then the claim holds since $\mid f \mid :I_{\mu _{\bar\nu}} \rightarrow I_{\mu _\nu}$ is $\Sigma _1$-elementary. If $\mu _\tau =\mu _\nu$ and $n(\tau )=n(\nu )$, then $P_\tau \subseteq P_\nu$. I.e. $\tau$ is dependent on $\nu$. Thus $\bar \tau$ is not independent. So let $\mu := \mu _\tau = \mu _\nu$, $n:=n(\tau )<n(\nu )$ and $\tau \in S^+ \cup  \widehat{S}$ be independent. Then, by the definition of the parameters, $\alpha _\tau$ is the $n$-th projectum of $\mu$.
\smallskip \\
Let 
$$\gamma _\beta :=crit(f_{(\beta ,0,\tau )})< \alpha _\tau$$
for a $\beta$ and
\smallskip 

$H_\beta$ $:=$ the $\Sigma _n$-hull of $\beta \cup P_\tau \cup \{ \alpha ^\ast _\mu , \tau \}$ in $I_\mu$.
\smallskip \\
I.e. $H_\beta = h^n_\mu [\omega \times (J^X_\beta \times \{ \alpha ^\prime _\mu ,\tau ^\prime ,P^\prime _\tau \} )]$ where
\smallskip 

$\alpha ^\prime _\mu$ $:=$ minimal such that $h^n_\mu (i,\alpha ^\prime _\mu )=\alpha ^\ast _\mu$ for an $i \in \omega$
\smallskip

$P^\prime _\tau$ $:=$ minimal such that $h^n_\mu (i,P^\prime _\tau )=P_\tau$ for an $i \in \omega$

$\tau ^\prime$ $:=$ minimal such that $h^n_\mu (i,\tau ^\prime )=\tau$ for an $i \in \omega$
(rsp. $\tau ^\prime := 0$ for $\tau =\mu$).
\smallskip \\
For the standard parameters are in $P_\tau$. 
\smallskip \\
so $H_\beta$ is $\Sigma _n$-definable over $I_\mu$ with the parameters $\{ \beta ,\tau ,\alpha ^\ast _\mu \} \cup P_\tau$ . Let
\smallskip 

$\rho$ $:=$ $\alpha _\tau$ $=$ the $n$-th projectum of $\mu$
\smallskip 

$A$ $:=$ the $n$-th standard code of $\mu$
\smallskip 

$p$ $:=$ $\langle \alpha ^\prime _\mu ,\tau ^\prime ,P^\prime _\tau \rangle$.
\smallskip \\
So $H_\beta \cap J^X_\rho$ is $\Sigma _0$-definable over $\langle I^0_\rho ,A \rangle$ with parameters $\beta$ and $p$. (fine structure theory!)
\smallskip \\
And $\gamma _\beta$ is defined by
$$\gamma _\beta \not \in H_\beta \quad and \quad (\forall \delta \in \gamma _\beta )(\delta \in H_\beta ).$$
I.e.  $\gamma _\beta$ is also $\Sigma _0$-definable over $\langle I^0_\rho ,A \rangle$ with parameters $\beta$ and $p$.
\smallskip \\
Let $f_0:=f_{(\beta ,0,\tau )}$ for a $\beta$, $\bar \tau _0 := d(f_0)<\alpha _\tau$ and $\gamma := crit(f_0) < \alpha _\tau$. Let $f_1:=f_{(\beta ,\gamma ,\tau )}$, $\bar \tau _1 := d(f_1)<\alpha _\tau$ and $\delta := crit(f_1) < \alpha _\tau$. Then $\mu _{\bar\tau _1}$ is the direct successor of $\mu _{\bar\tau _0}$ in $K_\tau$. So $f_{(\beta ,\gamma ,\bar \tau _1 )}=id_{\bar \tau _1}$. Hence $\mu _\eta =\mu _{\bar\tau _1}$ holds for the minimal $\eta \in S^+ \cup S^0$ such that $\gamma < \eta \sqsubseteq \delta$. Thus
\smallskip

$$\mu ^\prime \in K^+_\tau := K_\tau -(Lim(K_\tau ) \cup \{ min(K_\tau )\} )$$
$$\Leftrightarrow$$

$(\exists \beta ,\gamma ,\delta ,\eta )(\gamma =\gamma _\beta$ $and$ $\delta =\gamma _{(\gamma _\beta +1)}$
\smallskip

$and$ $\eta \in S^+ \cup S^0$ minimal such that $\gamma < \eta \sqsubseteq \delta$
$and$ $\mu ^\prime = \mu _\eta )$
\smallskip \\
Therefore, $K^+_\tau$ is $\Sigma _1$-definable over $\langle I^0_\rho ,A \rangle$ with parameter $p$.
\smallskip \\
Now, consider $\langle I^0_{\alpha _\tau},K_\tau  \rangle \models \varphi (x)$ where $\varphi$ is a $\Sigma _1$ formula. Then, since $K_\tau $ is unbounded in $\alpha _\tau$,
\smallskip
\begin{center}
$\langle I^0_{\alpha _\tau},K_\tau  \rangle \models \varphi (x)$
$$\Leftrightarrow$$
$(\exists \gamma )(\gamma \in K^+_\tau$ $and$ $\langle I^0_{\alpha _\gamma},K_\gamma  \rangle \models \varphi (x))$.
\end{center}
So $\langle I^0_{\alpha _\tau},K_\tau \rangle \models \varphi (x)$ is $\Sigma _1$ over $\langle I^0_\rho ,A \rangle$ with parameter $p$, rsp. $\Sigma _{n+1}$ over $I_\mu$ with parameters $\alpha ^\ast _\mu ,\tau ,P_\tau$. But since $n=n(\tau )<n(\nu )$, $f$ is at least $\Sigma _{n+1}$-elementary. In addition $f(\alpha ^\ast _{\bar \tau} )=\alpha ^\ast _\tau$, $f(\bar \tau )=\tau$, $f(P_{\bar \tau} )=P_\tau$. So, for $x \in rng(f)$, $\langle I^0_{\alpha _{\bar \tau}},K_{\bar \tau} \rangle \models \varphi (f^{-1}(x))$ holds iff $\langle I^0_{\alpha _\tau},K_\tau \rangle \models \varphi (x)$. 
$\Box$
\medskip \\
{\bf Theorem 29}
\smallskip \\
$\frak{M} := \langle S,\vartriangleleft ,\frak{F},D \rangle$ is a $\kappa$-standard morass.
\smallskip \\
{\bf Proof:}
Set 
$$\sigma _{(\xi ,\nu )} (i)=h^{n(\nu )}_\nu (i,\langle \xi ,\alpha _\nu ^\ast ,p_\nu  \rangle ).$$
Then $D$ is uniquely determined by the axioms of standard morasses and
\smallskip 

(1) $D^\nu$ is uniformly definable over $\langle J^X_\nu , X \upharpoonright \nu ,X_\nu \rangle$
\smallskip 

(2) $X_\nu$ is uniformly definable over $\langle J^D_\nu ,D_\nu ,D^\nu \rangle$.
\smallskip \\
(1) is clear. For (2), assume first that $\nu \in  \widehat{S}$ and $f_{(0,q_\nu ,\nu )}=id_\nu$. Since the set $\{ i \mid \sigma _{(q_\nu ,\nu )}(i) \in X_\nu \}$ is $\Sigma _{n(\nu )}$-definable over $\langle J^X_\nu ,X \upharpoonright \nu ,X_\nu \rangle$ with the parameters $p_\nu , \alpha _\nu ^\ast ,q_\nu$, there is a $j \in \omega$ such that
$$\sigma _{(q_\nu ,\nu )}(\langle i,j \rangle ) \quad \hbox{existiert} \Leftrightarrow \sigma _{(q_\nu ,\nu )}(i) \in X_\nu .$$
Using this $j$, we have
$$X_\nu = \{ \sigma _{(q_\nu ,\nu )}(i)\mid \langle i,j \rangle \in dom(\sigma _{(q_\nu ,\nu )}) \} .$$
So, in case that $f_{(0,q_\nu ,\nu )}=id_\nu$, there is the desired definition of $X_\nu$.
\smallskip \\
Let $\nu \in  \widehat{S}$, $f_{(0,q_\nu ,\nu )}:\bar\nu \Rightarrow \nu$ cofinal and $f(\bar q)=q_\nu$. Then $f_{(0, \bar q ,\bar \nu )}=id_{\bar\nu}$. And by lemma 6 (b) of [Irr2], $\bar q = q_{\bar \nu}$. So, if $\bar \nu =\nu$, then $f_{(0,q_\nu ,\nu )}=id_\nu$. Thus let $\bar \nu < \nu$. Then $f_{(0,q_\nu ,\nu )}(x)=y$ is defined by: There is a $\bar \nu \leq \nu$ such that, for all $r,s \in \omega$,
$$\sigma _{(q_{\bar \nu} ,\bar \nu )}(r) \leq \sigma _{(q_{\bar \nu} ,\bar \nu )}(s) \Leftrightarrow \sigma _{(q_\nu ,\nu )}(r) \leq \sigma _{(q_\nu ,\nu )}(s)$$
holds $and$ for all $z \in J^X_{\bar \nu}$ there is an $s \in \omega$ such that
$$z=\sigma _{(q_{\bar \nu} ,\bar \nu )}(s)$$
$and$ there is an $s \in \omega$ such that
$$\sigma _{(q_{\bar\nu} ,\bar \nu )}(s)=x \Leftrightarrow \sigma _{(q_\nu ,\nu )}(s)=y$$.
\smallskip \\
And since $\langle J^X_\nu ,X_\nu \rangle$ is rudimentary closed,
$$X_\nu = \bigcup \{ f(X_{\bar \nu} \cap \eta ) \mid \eta < \bar \nu \} .$$ 
Finally, if $\nu \in  \widehat{S}$ and $f_{(0,q_\nu ,\nu )}$ is not cofinal in $\nu$, then $C_\nu$ is unbounded in $\nu$ and
$$X_\nu = \bigcup \{ X_\lambda \mid \lambda \in C_\nu \} $$
by the coherence of $L_\kappa [X]$.
\smallskip \\
So (2) holds. From this, (DF)$^+$ follows.
\smallskip \\
By (1) and (2), $J^X_\nu =J^D_\nu$ for all $\nu \in Lim$, and for all $H \subseteq J^X_\nu =J^D_\nu$
$$H \prec _1 \langle J^X_\nu ,X \upharpoonright \nu \rangle \Leftrightarrow H \prec _1 \langle J^D_\nu ,D_\nu \rangle .$$
Now, we check the axioms.
\smallskip\\
(MP) and (MP)$^+$
\smallskip\\
$\mid f_{(0,\xi ,\nu )} \mid$ is the uncollapse of $h^{n[\nu )}_{\mu _\nu}[\omega \times \{ \xi ^\ast ,\nu ^\ast ,\alpha ^\ast _\nu ,\alpha ^{\ast\ast}_{\mu _\nu},P^\ast _\nu \} ^{<\omega}]$ where $\xi ^\ast$ is minimal such that $h^{n (\nu )-1}_{\mu _\nu}(i, \xi ^\ast )=\xi$. Therefore, (MP) and (MP)$^+$ hold.
\smallskip\\
(LP1)
\smallskip\\
holds by (2) above.
\smallskip\\
(LP2)
\smallskip\\
This is lemma 26.
\smallskip \\
(CP1) and (CP1)$^+$
\smallskip\\
This follows from lemma 24 and the definition of $\sigma _{(\xi ,\nu )}$.
\smallskip\\
(CP2)
\smallskip\\
This is lemma 27.
\smallskip\\
(CP3) and (CP3)$^+$
\smallskip \\
Let $x \in J^X_\nu$, $i \in \omega$ and $y=h_{\nu ,B_\nu}(i,x)$. Since $C_\nu$ is unbounded in $\nu$, there is a $\lambda \in C_\nu$ such that $x,y \in J^X_\lambda$. By lemma 25, $B_\lambda = B_\nu \cap J^X_\lambda$. So $y=h_{\lambda ,B_\lambda}(i,x)$.
\smallskip \\
(DP1)
\smallskip\\
holds by the definition of $\mu _\nu$.
\smallskip\\
(DF)
\smallskip \\
Let $\mu := \mu _\nu$, $k:=n(\mu )$ and
\smallskip

$\pi (n,\beta ,\xi )$ $:=$ the uncollapse of $h^{k+n}_\mu [\omega \times (J^X_\beta \times \{ \alpha _\mu ^{\ast\ast} ,p^\ast _\mu , \xi ^\ast \}^{<\omega} )]$
\smallskip\\
where 
\smallskip

$\xi ^\ast :=$ minimal such that $h^{k+n-1}_\mu (i,\xi ^\ast )=\xi$ for an $i \in \omega$
\smallskip

$p^\ast _\mu :=$ minimal such that $h^{k+n-1}_\mu (i,p _\mu ^\ast )=p _\mu$ for some $i \in \omega$
\smallskip

$\alpha _\mu ^{\ast\ast} :=$ minimal such that $h^{k+n-1}_\mu (i,\alpha _\mu ^{\ast\ast} )=\alpha _\mu ^\ast$ for some $i \in \omega$.
\smallskip \\
Prove 
$$\mid f^{1+n}_{(\beta ,\xi ,\mu )} \mid = \pi (n,\beta ,\xi ).$$
for all $n \in \omega$ by induction.
\smallskip\\
For $n=0$, this holds by definition of $f^1_{(\beta , \xi ,\mu )}=f_{(\beta , \xi ,\mu )}$.So assume that $\mid f^m_{(\beta ,\xi ,\mu )} \mid = \pi (m-1,\beta ,\xi )$ is already proved for all $1\leq m \leq n$. Then, by definition of $\tau (m,\mu )$,
\smallskip 

$\alpha _{\tau (m,\mu )}$ $=$ the $(k+m-1)$-th projectum of $\mu$.
\smallskip \\
Let $\pi (n,\beta ,\xi ):I_{\bar \mu} \rightarrow I_\mu$. Then
\smallskip

$(\ast ) \quad$ $\xi (m,\mu )=\pi(n,\beta ,\xi )\xi (m,\bar \mu )$ for all $1 \leq m \leq n$:
\smallskip \\
Let $\pi := \pi (n,\beta ,\xi )$, $\alpha := \pi ^{-1}[\alpha _{\tau (m,\mu )} \cap rng(\pi )]$, $\rho := \pi (\alpha )$
\smallskip 

$r:=$ minimal such that $h^{k+m-2}_\mu (i,r)=p_\mu$ for an $i \in \omega$
\smallskip 

$\alpha ^\prime :=$ minimal such that $h^{k+m-2}_\mu (i,\alpha ^\prime )=\alpha ^\ast _\mu$ for an $i \in \omega$
\smallskip 

$p:=$ the $(k+m-1)$-th parameter of $\mu$
\smallskip \\
and $\pi (\bar r)=r$, $\pi (\bar p)=p$, $\pi (\bar \alpha ^\prime )=\alpha ^\prime$. Let $\bar \xi := \xi (m,\bar \mu )$. Then $\bar p =h^{k+m-1}_{\bar \mu}(i,\langle \bar x , \bar \xi ,\bar r,\bar \alpha ^\prime \rangle )$ for a $\bar x \in J^X_\alpha$, because $\alpha = \alpha _{\tau (m,\bar \mu )}$. So $p=h^{k+m-1}_\mu (i,\langle x , \xi , r, \alpha ^\prime \rangle )$ where $\pi (\bar x)=x$ and $\pi (\bar \xi )=\xi$. Thus $h^{k+m-1}_\mu [\omega \times (J^X_{\alpha _{\tau (m,\mu )}} \times \{ \alpha ^\prime ,r,\xi \} ^{<\omega})]=J^X_\mu$ by definition of $p$. So $\xi (m,\mu )\leq \xi$. Assume $\xi (m,\mu )<\xi$. Then $I_\mu \models (\exists \eta < \xi )(\exists i \in \omega )(\exists x \in J^X_\rho )(\xi = h^{k+m-1}_\mu (i, \langle x , \eta , r, \alpha ^\prime \rangle )$. So $I_{\bar \mu}\models (\exists \eta < \bar \xi )(\exists i \in \omega )(\exists x \in J^X_\alpha )(\bar\xi = h^{k+m-1}_{\bar\mu} (i, \langle x , \eta , \bar r, \bar\alpha ^\prime \rangle )$. But this contradicts the definition of $\bar \xi =\xi (m,\bar \mu )$.
\smallskip \\
So, for all $1 \leq m \leq n$,
$$\xi (m,\mu ) \in rng(\pi (n,\beta ,\xi )).$$
In addition, for all $\beta < \alpha _{\tau (m,\mu )}$,
$$d(f^m_{(\beta ,\xi (m,\mu ),\mu )})<\alpha _{\tau (m,\mu )} .$$
Consider $\pi := \pi (m-1,\beta ,\xi ) = \mid f^m_{(\beta ,\xi ,\mu )} \mid$ where $\xi =\xi (m,\mu )$. Then $\pi :I_{\bar \mu} \rightarrow I_\mu$ is the uncollapse of $h^{k+m-1}_\mu [\omega \times (\beta \times \{ \xi ,\alpha ^\prime ,r \} ^{<\omega})]$ where 
\smallskip 

$r:=$ minimal such that $h^{k+m-2}_\mu (i,r)=p_\mu$ for some $i \in \omega$
\smallskip 

$\alpha ^\prime :=$ minimal such that $h^{k+m-2}_\mu (i, \alpha ^\prime )= \alpha ^\ast _\mu$ for some $i \in \omega$.
\smallskip \\
And $h^{k+m-1}_{\bar\mu} [\omega \times (\beta \times \{ \bar\xi , \bar\alpha ^\prime ,\bar r \} ^{<\omega})]=J^X_{\bar \mu}$ where $\pi (\bar \xi )=\xi$, $\pi (\bar \alpha ^\prime )=\alpha ^\prime$ and $\pi (\bar r)=r$. Assume $\alpha _{\tau (m,\mu )} \leq \bar \mu < \mu$. Then there were a function over $I_{\bar \mu}$ from $\beta < \alpha _{\tau (m,\mu )}$ onto $\alpha _{\tau (m,\mu )}$. This contradicts the fact that $\alpha _{\tau (m,\mu )}$ is a cardinal in $I_\mu$. If $\bar \mu = \mu$, then $f^m_{(\beta ,\bar \xi ,\mu )}=id_\mu$. This contadicts the minimality of $\tau (m,\mu )$.  
\smallskip \\
Since $\xi (m,\mu ) \in rng(\pi (n,\beta ,\xi ))$, we can prove 
$$rng(\pi (n,\beta ,\xi )) \cap J^D_{\alpha _{\tau (m,\mu )}} \prec _1 \langle J^D_{\alpha _{\tau (m,\mu )}},D_{\alpha _{\tau (m,\mu )}},K^m_\mu \rangle $$
for all $1 \leq m \leq n$ as in lemma 28.
\medskip \\
We still must prove minimality.. Let $f \Rightarrow \mu$ and $\beta \cup \{ \xi \} \subseteq rng(f)$ such that
$$rng(f) \cap J^D_{\alpha _{\tau (m,\mu )}} \prec _1 \langle J^D_{\alpha _{\tau (m,\mu )}},D_{\alpha _{\tau (m,\mu )}},K^m_\mu \rangle $$
$$\xi (m,\mu )\in rng(f)$$
holds for all $1 \leq m \leq n$. Show that $f$ is $\Sigma _{k+n}$-elementary and that the first standard parameters including the $(k+n-1)$-th are in $rng(f)$. That suffices because $\pi (n,\beta ,\xi )$ is minimal.
\smallskip \\
Let $p^{k+m}_\mu$ be the $(k+m)$-th standard parameter of $\mu$.
\smallskip \\
Prove, by induction on $0 \leq m \leq n$,
\smallskip 

$f$ is $\Sigma _{k+m}$-elementary
\smallskip

$p^1_\mu ,\dots ,p^{k+m-1}_\mu \in rng(f)$.
\smallskip \\
For $m=0$, this is clear because $f \Rightarrow \mu$. So assume it to be proved for $m<n$ already. Then let $\alpha := \alpha _{\tau (m+1,\mu )}$ and $\bar \alpha = f^{-1}[\alpha \cap rng(f)]$. Consider $\pi := (f \upharpoonright J^D_{\bar \alpha}): \langle J^D_{\bar \alpha},D_{\bar \alpha},\bar K \rangle \rightarrow \langle J^D_\alpha ,D_\alpha ,K^{m+1}_\mu \rangle$. Construct a  $\Sigma _{k+m+1}$-elementary extension $\tilde \pi$ of $\pi$. To do so, set
$$f_\beta =f^{m+1}_{(\beta , \xi (m+1,\mu ),\mu )}$$
$$\mu (\beta )=d(f_\beta )$$
$$H= \bigcup \{ f_\beta [rng(\pi ) \cap J^D_{\mu (\beta )}]\mid \beta < \alpha \} .$$  
Then $H \cap J^D_\alpha =rng(\pi )$. For, $rng(\pi ) \subseteq H \cap J^D_\alpha$ is clear because $f_\beta \upharpoonright J^D_\beta = id \upharpoonright J^D_\beta$. So let $y \in H \cap J^D_\alpha$. I.e. $y=f_\beta (x)$ for some $x \in rng(\pi )$ and a $\beta < \alpha$. Let 
$K^+ =K^{m+1}_\mu -Lim(K^{m+1}_\mu )$ and $\beta (\eta )= sup \{ \beta \mid f^{m+1}_{(\beta ,\xi (m+1,\eta ), \eta )}\neq id_\eta \}$. Then
$$\langle J^D_\alpha , D_\alpha ,K^{m+1}_\mu \rangle \models (\exists y)(\exists \eta \in K^+ )(y=f^{m+1}_{(\beta ,\xi (m+1,\eta ), \eta )}(x) \in J^D_{\beta (\eta )}).$$
Since $rng(\pi ) \prec _1 \langle J^D_\alpha ,D_\alpha ,K^{m+1}_\mu \rangle$, $y=f^{m+1}_{(\beta ,\xi (m+1,\eta ), \eta )}(x) \in rng(\pi )$ if $x \in rng(\pi )$ for such an $\eta$. But since $y=f^{m+1}_{(\beta ,\xi (m+1,\eta ), \eta )}(x) \in J^D_{\beta (\eta )}$, we get $f_\beta (x) = f^{m+1}_{(\beta ,\xi (m+1,\eta ), \eta )}(x) \in rng(\pi )$.
\smallskip \\
Show $H \prec _{k+m+1} I_\mu$. Since $f^{m+1}_{(\beta , \xi ,\mu )}=\pi (m,\beta ,\xi )$, $\alpha _{\tau (m+1,\mu )}$ is the $(k+m)$-th projectum of $\mu$. Like in $(\ast )$ above, we can show that the $(k+m)$-th standard parameter $p^{k+m}_\mu$ of $\mu$ is in $rng(f_\beta )$. Now, let $I_\mu \models (\exists x)\varphi (x,y,p^1_\mu ,\dots ,p^{k+m}_\mu )$ where $\varphi$ is a $\Pi _{k+m}$ formula and $y \in H \cap J^D_\alpha$. Since $f_\beta$ is $\Sigma _{k+m}$-elementary, the following holds:
$$ I_\mu \models (\exists x)\varphi (x,y,p^1_\mu ,\dots ,p^{k+m}_\mu ) \Leftrightarrow (\exists \gamma \in K^{m+1}_\mu )(\exists x)(I_\gamma \models \varphi (x,y,p^1_\gamma ,\dots ,p^{k+m}_\gamma )).$$ 
And since $rng(\pi ) \prec _1 \langle J^D_\alpha ,D_\alpha ,K^{m+1}_\mu \rangle$, 
$$rng(\pi ) \models (\exists \gamma \in K^{m+1}_\mu )(\exists x)(I_\gamma \models \varphi (x,y,p^1_\gamma ,\dots ,p^{k+m}_\gamma )).$$
Thus there is such an $x$ in $rng(\pi )$ and therefore in $H$.
\smallskip \\
Let $\tilde \pi$ be the uncollapse of $H$. Then $\tilde \pi$ is $\Sigma _{k+m}$-elementary and, since $p^1_\mu ,\dots ,p^{k+m}_\mu \in rng(f_\beta )$ for all $\beta < \alpha$, we have $p^1_\mu ,\dots ,p^{k+m}_\mu \in rng(\pi )=H$. In addition, by the induction hypothesis, $f$ is $\Sigma _{k+m}$-elementary and $p^1_\mu ,\dots ,p^{k+m-1}_\mu \in rng(f)$. Again as in $(\ast )$ above, we can show that $p^{k+m}_\mu \in rng(f)$ using  $\xi (m+1,\mu )\in rng(f)$. But since $\tilde \pi$ and $f$ are the same on the $(k+m)$-th projectum, we get $\tilde \pi =f$.
\smallskip \\
(SP) follows from $\mid f^{1+n}_{(\beta ,\xi ,\mu )} \mid = \pi (n,\beta ,\xi )$, because for all $\nu \sqsubset \tau \sqsubseteq \mu _\nu$ such that $\tau \in S^+$ (rsp. $\tau =\nu$) the following holds:
$$p_\tau \in rng(\pi (n,\beta ,\xi )) \Leftrightarrow \xi _\tau \in rng(\pi (n,\beta ,\xi )).$$
This may again be shown as $(\ast )$.
\smallskip \\
(DP2)
\smallskip\\
is like $(\ast )$ in (DF).
\smallskip\\
(DP3)
\smallskip \\
(a) is clear.
\smallskip \\
(b) was already proved with (DF)$^+$.
\smallskip\\
$\Box$
\bigskip \\
{\bf Theorem 30}
\smallskip \\
Let $\langle X_\nu \mid \nu \in S^X \rangle$ be such that
\smallskip

(1) $L[X] \models S^X =\{ \beta (\nu )\mid \nu$ singular$\}$
\smallskip 

(2) $L[X]$ is amenable
\smallskip 

(3) $L[X]$ has condensation
\smallskip

(4) $L[X]$ has coherence.
\smallskip\\
Then there is a sequence $C = \langle C_\nu \mid \nu \in  \widehat{S} \rangle$ such that
\smallskip 

(1) $L[C]=L[X]$
\smallskip

(2) $L[C]$ has condensation
\smallskip

(3) $C_\nu$ is club in $J^C_\nu$ w.r.t. the canonical well-ordering $<_\nu$ of $J^C_\nu$
\smallskip 

(4) $otp(\langle C_\nu ,<_\nu \rangle )>\omega \Rightarrow C_\nu \subseteq \nu$
\smallskip

(5) $\mu \in Lim(C_\nu )$ $\Rightarrow$ $C_\mu =C_\nu \cap \mu$,
\smallskip

(6) $otp(C_\nu )<\nu$.
\smallskip \\
{\bf Proof:} First, construct from $L[X]$ a standard morass as in theorem 29. Then construct a inner model $L[C]$ from it as in [Irr2]. $\Box$
\section*{References}
[BJW] A. Beller, R.-B. Jensen, P. Welch: {\bf Coding the Universe}, London Math\-ematical Society Lecture Notes Series, Cambridge University Press, London, 1982
\smallskip \\
{[Dev1]} K. Devlin: {\bf Aspects of Constructibility}, Lecture Notes in Mathematics, no. 354, Springer-Verlag, Berlin, 1973
\smallskip \\
{[Dev2]} K. Devlin: {\bf Constructibility}, Springer-Verlag, Berlin, 1984
\smallskip \\
{[DJS]} H.-D. Donder, R.-B. Jensen, L.J. Stanley: Condensation-Coherent Global Square Systems, {\bf Recursion theory}, Proceedings of Symposia in Pure Mathematics, vol. 42, American Mathematical Society, Providence, 1985, 237--258
\smallskip \\
{[Irr1]}B. Irrgang: {\bf Kondensation und Moraste}, Dissertation, M\"unchen, 2002
\smallskip\\
{[Irr2]} B. Irrgang: {\bf Proposing $(\omega _1 ,\beta )$-Morasses, $\omega _1 \leq \beta$}
\smallskip\\
{[Jen1]} R.-B. Jensen: The Fine Structure of the Constructible Hierarchy, {\bf Annals of Mathematical Logic}, vol. 4 (1972), 229--308
\smallskip \\
{[Jen2]} R.-B. Jensen: {\bf Higher-Gap Morasses}, hand-written notes, 1972/73
\smallskip \\
\end{document}